\documentclass{amsart}
\usepackage{fullpage}
\usepackage{amsmath}
\usepackage{amssymb}
\usepackage{mathabx}
\usepackage{todonotes}
\usepackage{tikz-cd}

\begin{document}

\newcommand{\+}{\boxplus}
\newcommand{\x}{\cdot}
\newcommand{\del}{\setminus}
\newcommand{\ignore}[1]{}
\newtheorem{example}{Example}
\newtheorem{lemma}{Lemma}
\newtheorem{theorem}{Theorem}
\newtheorem{corollary}{Corollary}
\newtheorem{conjecture}{Conjecture}
\newtheorem{claim}{Claim}[theorem]
\newcommand{\Z}{\mathbb{Z}}
\newcommand{\N}{\mathbb{N}}
\newcommand{\one}{{\bf 1}}
\newcommand{\R}{\mathbb{R}} 

\title{Field extensions, Derivations, and \\
Matroids over skew hyperfields}
\begin{abstract} 
We show that a field extension $K\subseteq L$ in positive characteristic $p$ and elements $x_e\in L$ for $e\in E$  gives rise to a matroid $M^\sigma$ on ground set $E$ with coefficients in a certain skew hyperfield $L^\sigma$. This skew hyperfield $L^\sigma$ is defined in terms of $L$ and its Frobenius action $\sigma:x\mapsto x^p$.  The matroid underlying $M^\sigma$ describes the algebraic dependencies over $K$ among the $x_e\in L$ , and $M^\sigma$ itself comprises, for each $m\in \Z^E$, the linear space of $K$-derivations of $K\left(x_e^{p^{m_e}}: e\in E\right)$.

The  theory of matroid representation over hyperfields was developed by Baker and Bowler for commutative hyperfields. We partially extend their theory to skew hyperfields. To prove the duality theorems we need,  we use a new axiom scheme in terms of {\em quasi-Pl\"ucker coordinates}. 
\end{abstract}

\author{Rudi Pendavingh}

\address[Rudi Pendavingh]{Department of Mathematics and Computer Science,
Eindhoven University of Technology, P.O.~Box 513, 5600 MB Eindhoven, The Netherlands}
\email{r.a.pendavingh@tue.nl}

\maketitle
\section{Introduction}
Let $K$ be a field and let $(x_e: e\in E)$ be elements from an extension field $L$. A subcollection $(x_e: e\in F)$ is {\em algebraically dependent over $K$} if there is a polynomial $q\in K[X_e: e\in F]$ so that $q(x_e: e\in F)=0$. By a theorem of Steinitz, the set 
$\mathcal{I}:=\{F\subseteq E: (x_e: e\in F)\text{ algebraically independent over }K\}$
satisfies
\begin{itemize}
\item[(I0)] $\emptyset\in \mathcal{I}$
\item[(I1)] if $A\in \mathcal{I}$ and $B\subseteq A$, then $B\in \mathcal{I}$
\item[(I2)] if $A, B\in\mathcal{I}$ and $|A|<|B|$, then $A\cup\{e\}\in\mathcal{I}$ for some $e\in B\del A$
\end{itemize}
Algebraic independence has these properties in common with linear independence. This formalizes the analogy beween algebraic closure and linear span, transcendence degree and dimension of a linear space, and in general gives a geometric perspective on field extensions. 

A {\em matroid} is a pair $M=(E,\mathcal{I})$ where $E$ is a finite set and $\mathcal{I}$ is any set of subsets of $E$ satisfying (I0), (I1), and (I2). The above pair $K, x$ thus gives an {\em algebraic matroid} $M(K,x)$, and a collection of vectors $(v_e)_{e\in E}$ will determine a {\em linear matroid} on $E$.

Taking poetic license, a matroid may be described as a {\em linear space without coefficients}. In a linear space over a field $K$ of dimension $d$, any subset of $d$ vectors is associated with a value in $K$, the determinant. The corresponding matroid merely distinguishes between bases and non-bases. There are less Spartan matroid variants, such as oriented matroids and valuated matroids, which can be seen as matroids with coefficients in the set of signs $\{+,-\}$ and in a linearly ordered group, respectively. This intuitive perspective was developed rigorously by Dress and Wenzel \cite{DW1991}, who defined {\em matroids with coefficients} from a {\em fuzzy field}, and more recently by Baker and Bowler \cite{BakerBowler2017}, who defined {\em matroids over hyperfields}. In both approaches,  linear spaces as well as oriented-, valuated-, and ordinary matroids are matroids with coefficients in a corresponding fuzzy field or hyperfield. 

Hyperfields generalize fields, and their more relaxed additive structure translates to a richer collection of homomorphisms. A hyperfield homomorphism $f:H\rightarrow H'$ induces a map $f_*$ which takes a matroid over $H$ and turns it into a matroid over $H'$, simply by applying $f$ to the coefficients. This elegantly describes how a matroid with coefficients in a field $K$ (essentially a collection of vectors in a $K$-vector space) gives rise to an ordinary matroid. From any field $K$, there is a hyperfield homomorphism $\kappa$ to the {\em Krasner hyperfield} $\mathbb{K}=\{0,1\}$, which maps $0\mapsto 0$ and all nonzero $x\mapsto 1$. The induced  map $\kappa_*$ is a forgetful operation which retains only the destinction between bases and nonbases, independent and dependent sets. If the ordinary matroid that arises from applying $\kappa_*$ appears too coarse an abstraction, then one may consider a hyperfield homomorphism from $K$ to a more 
detailed hyperfield. For example, the natural homomorphism from the reals to the hyperfield of signs induces the map from collections of vectors in Euclidean space to oriented matroids.

So in the study of linear spaces, matroids over hyperfields may serve to attain the `right' abstraction level. In relation to field extensions, they have a different role.  Unlike a linear space, a field extension is not itself a matroid over some hyperfield. The algebraic matroid $M(K,x)$ is, but the information on $(K,x)$ it contains is quite sparse. To illustrate, we know of no general method to decide if $N=M(K,x)$ for some $(K,x)$, given a matroid $N$. 

In \cite{BDP2018}, it was show that a pair $(K,x)$ also determines a matroid valuation of  $M(K,x)$, the {\em Lindstr\"om valuation}. That is, $(K,x)$ determines a matroid over the hyperfield $\Z_{\min}$ with underlying matroid $M(K,x)$. In this paper, we show that $M(K,x)$ can even be decorated with coefficients in a certain hyperfield $L^\sigma$, which is defined in terms of $L$ and the Frobenius action $\sigma: x\mapsto x^p$. The {\em left $L^\sigma$-matroid $M^\sigma(K,x)$} that arises is still a geometric object, but comprises more detailed information about the pair $(K,x)$, such as the space of $K$-derivations of $K(x_e:  e\in E)$.

If $K'$ is an extension field of $K$, then a {\em $K$-derivation of $K'$} is any map $D:K'\rightarrow K'$ which is trivial on $K$, is additive, and satisfies the {\em Leibnitz rule} $D(xy)=D(x)y+xD(y)$.  The collection $Der(K, K')$ of all $K$-derivations of $K'$ is a linear space whose dimension in equals the transcendence degree of $K'$ over $K$. If $K':=K(x_e:  e\in E)$, then the dimension of $Der(K, K')$ equals the rank of the algebraic matroid  $M(K,x)$. The linear space $Der(K, K')$ then induces a linear matroid $M'(K,x)$ on $E$ of the same rank as  $M(K,x)$, in which a set $B\subseteq E$ is a basis if and only if for each $u\in (K')^B$, there is a unique $K$-derivation $D$ of $K'$ such that $D(x_e)=u_e$ for all $e\in B$. Such a basis of $M$ is necessarily a basis of $M(K,x)$, but the converse need not be true. In other words, the matroid of derivations $M'(K,x)$ is a {\em weak image} of $M(K,x)$. 
 
For any $m\in\Z^E$ and $x\in L^E$, let $\sigma^m(x):=(\sigma^{m_e}(x_e): e\in E)$. Passing from $x$ to $\sigma^m(x)$ does not affect algebraic dependence, and we have $M(K,x)=M(K,\sigma^m(x))$ for any $m\in \Z^E$. The matroid $M^\sigma(K,\sigma^m(x))$ arises from $M^\sigma(K,x)$ by {\em rescaling}, an operation which is defined generally for matroids over hyperfields. The matroid of derivations $M'(K,\sigma^m(x))$ in general does not equal $M'(K,x)$, and there is no easy relation between the two. But $M^\sigma(K, \sigma^m(x))$ does determine the space of $K$-derivations of $K(\sigma^m(x))$, and hence via rescaling, $M^\sigma(K,x)$ describes both this space of derivations and its matroid $M'(K,\sigma^m(x))$ for each $m\in \Z^E$. We have the following diagram. 
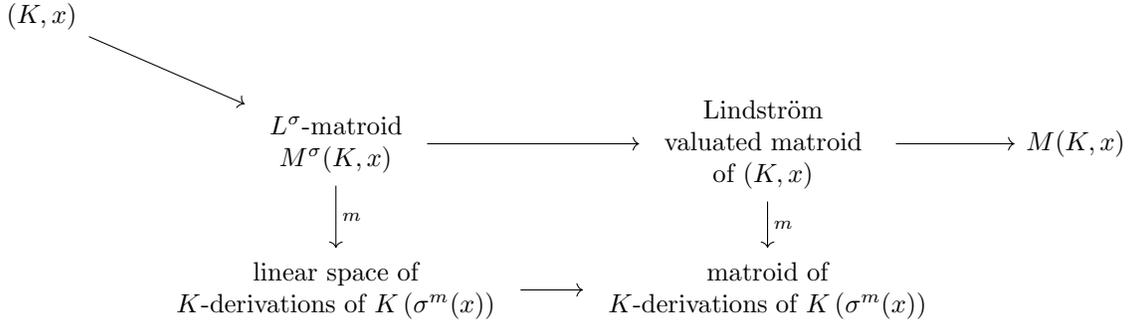
\begin{center}
\begin{tikzcd}
  (K, x)\arrow[rd]\\  
  
  & \text{\begin{tabular}{c}$L^\sigma$-matroid\\ $M^\sigma(K,x)$\end{tabular}} \arrow [r]\arrow[d,"m"]& \mbox{\begin{tabular}{c}Lindstr\"om \\ valuated matroid \\ of $(K,x)$\end{tabular} }\arrow[d, "m"]\arrow [r]& M(K,x)\\
  &\text{\begin{tabular}{c}linear space of \\$K$-derivations of $K\left(\sigma^m(x)\right)$\end{tabular}} \arrow[r] &\text{\begin{tabular}{c}matroid of \\$K$-derivations of $K\left(\sigma^m(x)\right)$\end{tabular}}\\
\end{tikzcd}
\end{center}
With the exception of  $(K,x)$ on the left, each node in this diagram is a matroid over a hyperfield, and each arrow represents a well-defined forgetful operation. Horizontal arrows indicate the application of a hyperfield homomorphism to the matroid coefficients, preserving the underlying matroid. Vertical arrows represent a new operation on matroids over certain hyperfields, which in general replaces the underlying matroid with a weak image of that matroid, and restricts the hyperfield to a sub-hyperfield.

As the diagram indicates, $M^\sigma(K,x)$ determines a map 
$$\mathcal{V}:m\mapsto  \left\{\text{$K$-derivations of $K\left(\sigma^{m}(x)\right)$}\right\}.$$
Essentially this object was called a {\em Frobenius flock} in \cite{BDP2018}. It was show in that paper that the related {\em matroid flock} $\mathcal{M}:m\mapsto M(\mathcal{V}_m)$ is a cryptomophic description of a matroid valuation of $M(K,x)$, which we named the Lindstr\"om valuation. This definition of the Lindstr\"om valuation via flocks was somewhat indirect, but shortly after a preprint of \cite{BDP2018} appeared on arXiv, Dustin Cartwright presented a direct construction of the Lindstr\"om valuation in \cite{Cartwright2017}. 

So matroid flocks are cryptomorphic to valuated matroids, and valuated matroids `are' matroids represented over the tropical hyperfield. Matroid flocks arise by a forgetful operation from Frobenius flocks. This suggested that perhaps, Frobenius flocks are also cryptomorphic to matroids represented over a certain hyperfield, and that the operation by which a Frobenius flock begets a matroid flock is just the pushing forward along an appropriate hyperfield homomorphism. In this paper, we show that this is exactly the case, the cryptomorphic description of the Frobenius flock of $(K, x)$ being the $L^\sigma$-matroid $M^\sigma(K,x)$. Rather than constructing $M^\sigma(K,x)$ via the Frobenius flock, we use the approach of Cartwright, and define  $M^\sigma(K,x)$ directly in terms of $(K,x)$.

The hyperfield $L^\sigma$ used to alternatively describe Frobenius flocks as left $L^\sigma$-matroids is not commutative. The theory of matroids over hyperfields was developed so far for commutative hyperfields. In the center of the theory of Baker and Bowler is the notion of a {\em Grassmann-Pl\"ucker function} of a matroid over a hyperfield, which generalizes the Pl\"ucker coordinates of a linear subspace. There is no proper analogue of the Grassmann-Pl\"ucker function in the context of skew hyperfields, just as there is no clean way to define the determinant of a matrix over a skew field.

However,  Gelfand, Gelfand, Retakh, and Wilson  \cite{GGRW2005} show that matrices over skew fields do admit {\em quasi-determinants}, which in the commutative setting equal ratios of certain adjacent determinants. Using this new concept, they also define {\em quasi-Pl\"ucker coordinates} for a matrix over a skew field, which are invariants of the linear space spanned by the rows of the matrix. As it turns out, this concept blends perfectly with matroids over hyperfields, and this allows us to  replace the Grassmann-Pl\"ucker functions with {\em quasi-Pl\"ucker coordinates} in the context of skew hyperfields. 

The structure of the paper is as follows. After giving preliminaries on matroids and hyperfields in Section \ref{sec:prelim}, we develop  {\em matroids over skew hyperfields} in Section \ref{sec:skew}. To demonstrate that the quasi-Pl\"ucker coordinates are natural in the context of matroids over skew hyperfields, we chose to make the presentation self-contained, but clearly most concepts and ideas in this section are adapted from Baker and Bowler \cite{BakerBowler2017}, Gelfand, Gelfand, Retakh, and Wilson  \cite{GGRW2005}, and others.
In Section \ref{sec:hyperfield}, we describe how to construct a {\em skew hyperfield of monomials} $H^\sigma$ from any hyperfield $H$ with automorphism $\sigma$. We describe the operation indicated by the vertical arrows in the diagram, which in general takes a matroid $M$ with coefficients in $H^\sigma$ and produces a matroid with coefficients in $H$, the {\em boundary matroid} $M_0$.
In Section \ref{sec:algebraic}, we show that each algebraic matroid representation $x$ in a field extension $L/K$ gives rise to a left $L^\sigma$-matroid, the {\em matroid of $\sigma$-derivations} $M^\sigma(K, x)$. The spaces of derivations as in the diagram arise from $M^\sigma(K,x)$ by rescaling and then taking the boundary matroid, so that $M^\sigma(K,x)$ determines the Frobenius flock. In general, a $H^\sigma$-matroid $M$ will determine a flock of $H$-matroids. In Section \ref{sec:flocks}, we prove that this flock in turn determines $M$.
In the final section of the paper, we make a few more related comments and present some conjectures.


\section{\label{sec:prelim}Preliminaries}
\subsection{Hypergroups, hyperrings, and hyperfields} 
A {\em hyperoperation} on $G$ is a map $\boxplus:G\times G\rightarrow 2^G$. Any hyperoperation induces a map $\overline{\boxplus}:2^G\times 2^G\rightarrow 2^G$ by setting 
$$X ~\overline{\boxplus} ~Y:=\bigcup\{x\boxplus y: x\in X, y\in Y\}.$$
Slightly abusing notation, one writes $x\+ Y:=\{x\}~\overline{\boxplus}~Y$, $X\+ y:=X~\overline{\boxplus}~\{y\}$, and $X\+Y:=X~\overline{\boxplus}~Y$.
The hyperoperation $\boxplus$ then is {\em associative}  if $x\+(y\+ z)=(x\+ y)\+z$ for all $x,y,z\in G$.

A {\em hypergroup} is a triple $(G, \boxplus, 0)$, where $0\in G$ and $\boxplus: G\times G\rightarrow 2^G\del\{\emptyset\}$ is an associative hyperoperation, such that
\begin{itemize}
\item[(H0)] $x\boxplus 0=\{x\}$
\item[(H1)] for each $x\in G$ there is a unique $y\in G$ so that $0\in x\boxplus y$. We write $-x:=y$
\item[(H2)] $x\in y\boxplus z$ if and only if $z\in x\boxplus (-y)$
\end{itemize}
If $G, H$ are hypergroups, then a map $f:G\rightarrow H$ is a {\em hypergroup homomorphism} if $f(x\+y)\subseteq f(x)\+f(y)$ for all $x,y\in G$, and $f(0)=0$. 

A {\em hyperring} is a tuple $(R, \cdot, \boxplus, 1, 0)$ so that
\begin{itemize}
\item[(R0)] $(R, \boxplus, 0)$ is a commutative hypergroup
\item[(R1)]  $(R^\star, \cdot, 1)$ is monoid, where we denote $R^\star:=R\setminus\{0\}$
\item[(R2)]  $0\cdot x=x\cdot 0=0$ for all $x\in R$
\item[(R3)]  $\alpha(x\boxplus y)=\alpha x\boxplus \alpha y$ and $(x\boxplus y)\alpha=x \alpha \boxplus y\alpha $ for all $\alpha,x,y\in R$
\end{itemize}
If $R, S$ are hyperrings, then $f:R\rightarrow S$ is a {\em hyperring homomorphism} if $f$ is a hypergroup homomorphism, $f(1)=1$, and $f(x\x y)=f(x)\x f(y)$ for all $x, y\in R$. 

A {\em skew hyperfield} is a hyperring such that $0\neq 1$, and each nonzero element has a multiplicative inverse. A {\em hyperfield} is then a skew hyperfield with commutative multiplication.  A (skew) hyperfield homomorphism is just a homomorphism of the underlying hyperrings.

The {\em Krasner hyperfield} is $\mathbb{K}=(\{0,1\}, \cdot, \boxplus, 1, 0)$, where $1\boxplus 1=\{0,1\}$. All hyperfields $H$ admit a hyperfield homomorphism $\kappa:H\rightarrow \mathbb{K}$ so that $\kappa(x)=1$ for all nonzero $x\in H$. Any skew field can be considered a skew hyperfield with hyperaddition $x\+y=\{x+y\}$. 

If  $(\Gamma, 0, +, <)$ is a linearly ordered abelian group, then $\Gamma_{\min}:=(\Gamma\cup\{\infty\}, 0,\infty, \odot, \boxplus)$ is a hyperfield, where we denoted $i\odot j:=i+j$ and 
$$i\boxplus j:=\left\{\begin{array}{ll} 
\{\min\{i,j\}\}&\text{if }i\neq j\\
\{m\in \Gamma: m\geq i\}\cup\{\infty\}& \text{if }i= j
\end{array}
\right.
$$
If $K$ is a (skew) field, then a map $\nu: K\rightarrow \Gamma\cup\{\infty\}$ is a {\em non-archimedean valuation} exactly if $\nu$ is a hyperfield homomorphism from $K$ to $\Gamma_{\min}$.
Replacing $\min$ with $\max$, $\geq$ with $\leq$, and $\infty$ with $-\infty$, we analogously obtain $\Gamma_{\max}$. In this paper, we use the hyperfield $\Z_{\min}$ as obtained from this construction. 

The smallest non-abelian group can be fitted with a hyperaddition to form a skew hyperfield. Consider $\mathbb{D}_3:=(D_3\cup \{0\}, \cdot, \boxplus, 1, 0)$, where 
$(D_3, \cdot, 1)$ is the dihedral group presented as $D_3=\{d_i: i\in \Z_6\}$ with $1:=d_0$, with multiplication and hyperaddition fixed by
$$d_i\cdot d_j =\left\{\begin{array}{ll} d_{i+j}&\text{if }i\in\{0,2,4\}\\ d_{i-j}&\text{if }i\in\{1,3,5\}\end{array}\right. \text{ and }
d_i\+ d_j =\left\{\begin{array}{ll} \{d_i\}&\text{if }j=i\\ \{d_i, d_j\}&\text{if }j=i+1   \\ \{d_i, d_{i+1},d_j\}&\text{if }j=i+2\\    D_3\cup \{0\}&\text{if }j=i+3  \end{array}\right.
$$
Verifying that $\mathbb{D}_3$ is indeed a skew hyperfield amounts to a finite check, which we omit.

\subsection{Matroids} A {\em matroid} is a pair $(E,\mathcal{C})$, where $E$ is a finite set and $\mathcal{C}$ is a set of subsets of $E$, such that 
\begin{itemize}
\item[(MC0)] $\emptyset\not\in \mathcal{C}$
\item[(MC1)] if $C, C' \in \mathcal{C}$ and $C\subseteq C'$, then $C=C'$
\item[(MC2)] for all  distinct $C, C'\in \mathcal{C}$ and all $e\in C\cap C'$, there exists a $C''\in \mathcal{C}$ such that $e\not \in C''\subseteq C\cup C'$
\end{itemize}
The elements of $\mathcal{C}$ are the {\em circuits} of the matroid $M=(E,\mathcal{C})$, and $E$ is the {\em ground set}. A subset $F$ of $E$ is {\em dependent} if $F\supseteq C$ for some $C\in \mathcal{C}$, and is {\em independent} otherwise. An inclusion-wise maximal independent set is called a {\em basis}. In a matroid $M$, all bases have the same cardinality, and this common cardinality is called the {\em rank} of $M$. 

In the context of a matroid $M$ with ground set $E$, we will write subsets of $E$ concisely as e.g. $Fabc:=F\cup \{a,b,c\}$. When we use this format, it is assumed implicitly that $a,b,c$ are distinct elements of $E\setminus F$. So a phrase `suppose $Fab$ is a basis of $M$' hides the more elaborate setup `suppose $F\subseteq E$, and $a,b$ are distinct elements of $E\setminus F$ so that $F\cup\{a,b\}$ is a basis of $M$'.

If $E$ is a finite set, $K$ is a field $V$ is a $K$-linear vector space, and $v_e\subseteq V$ for each $e\in E$, then for each $F\subseteq E$, the set $\{v_e: e\in F\}$ is either linearly dependent or independent over $K$. This distinction between dependent and independent sets is matroidal: if $\mathcal{C}$ denotes the set of inclusion-wise minimal nonempty sets $F$ corresponding to a dependent set of vectors $\{v_e: e\in F\}$, then $\mathcal{C}$ satisfies the circuit axioms (MC0), (MC1), and (MC2), and thus $M=(E, \mathcal{C})$ is a matroid. 

If $\mathcal{C}\subseteq 2^E$, then we say that a pair of distinct elements $C, C'\in \mathcal{C}$ is {\em modular} if $C\cup C'$ does not properly contain the union of two distinct elements of $\mathcal{C}$. Consider the {\em modular circuit elimination axiom}:
\begin{itemize}
\item[(MC2)'] for all  modular $C, C'\in \mathcal{C}$ and all $e\in C\cap C'$, there exists a $C''\in \mathcal{C}$ such that $e\not \in C''\subseteq C\cup C'$.
\end{itemize}
Then in the presence of (MC0) and (MC1), the ordinary circuit elimination axiom (MC2) is implied by its seemingly weaker modular counterpart (MC2)', so that we could alternatively define a matroid as a pair  $(E,\mathcal{C})$ for which (MC0), (MC1), and (MC2)' hold. The definition of {\em weak matroids over hyperfields} in \cite{BakerBowler2017} generalizes the {\em signed circuit axioms} for phased matroids given by Anderson and Delucchi \cite{AndersonDelucchi2012}, as well as the modular circuit axioms for ordinary matroids.

\ignore{
This above example may also be presented as follows. Let $L\subseteq K^E$ be a linear subspace. Call a subset $F\subseteq E$ {\em dependent} if there is a nonzero vector $x\in L$ so that $\underline{x}\subseteq F$, where $\underline{x}:=\{e\in E: x_e\neq 0\}$. It is straightforward that this notion of dependence is matroidal: if $x,y\in L$ are nonzero vectors so that $\underline{x}\neq \underline{y}$, and $e\in \underline{x}\cap\underline{y}$, then $z:= x_e^{-1}x-y_e^{-1}y\neq 0$ and $e\not \in \underline{z}\subseteq \underline{x}\cup\underline{y}$.
}


We refer to Oxley's book \cite{OxleyBook} for further matroid-related notation and results and to the paper of Baker and Bowler \cite{BakerBowler2017} for the theory of matroids over (commutative) hyperfields.


\section{\label{sec:skew} Matroids over skew hyperfields}
\subsection{Circuit axioms}
Let $H$ be a skew hyperfield, and let $E$ be a finite set. For any $X\in H^E$, let $\underline{X}:=\{e\in E: X_e\neq 0\}$ denote the {\em support} of $X$. A {\em  left $H$-matroid on $E$} is a pair $(E, \mathcal{C})$, where $\mathcal{C}\subseteq H^E$ satisfies the following {\em circuit axioms}.
\begin{itemize}
\item[(C0)] $0\not\in \mathcal{C}$.
\item[(C1)] if $X\in \mathcal{C}$ and $\alpha\in H^\star$, then $\alpha\x X\in \mathcal{C}$.
\item[(C2)] if $X,Y\in \mathcal{C}$ and $\underline{X}\subseteq \underline{Y}$, then there exists an $\alpha\in H^\star$ so that $Y = \alpha\x X$.
\item[(C3)] if $X,Y\in\mathcal{C}$ are a modular pair in $\mathcal{C}$ and $e\in E$ is such that $X_e=-Y_e\neq 0$, then there exists a $Z\in\mathcal{C}$ so that $Z_e=0$ and $Z\in X\+ Y$.
\end{itemize}
 In (C3), a pair  $X,Y\in\mathcal{C}$ is {\em modular} if $\underline{X},\underline{Y}$ are modular in $\underline{\mathcal{C}}:=\{\underline{X}: X\in \mathcal{C}\}.$ The notation $Z\in X\+ Y$ is short for $Z_f\in X_f\+ Y_f$ for all $f\in E$. 

A {\em  right $H$-matroid} is defined analogously, with $\alpha\x X$ replaced by $ X\x \alpha$ in (C1) and (C2). If $H$ is commutative, then left- and right $H$-matroids coincide, and we speak of {\em $H$-matroids}\footnote{In \cite{BakerBowler2017}, Baker and Bowler consider both {\em weak} and {\em strong} matroids over a hyperfield; our $H$-matroids are their {\em weak} $H$-matroids.}.

\ignore{
If $K$ is any skew field, then a linear subspace $L\subseteq K^E$ induces a left $K$-matroid $M(L)=(E,\mathcal{C})$ with circuits
$$\mathcal{C}:=\{X\in L\del\{0\}: \text{if }Y\in L\del\{0\}\text{ and }\underline{Y}\subseteq \underline{X}, \text{ then }Y=\alpha X\text{ for some }\alpha\in K\}$$ 
In turn, $L$ is determined by $M(L)$ as the linear subspace of $K^E$ spanned by $\mathcal{C}$.
}
Suppose $E$ is a finite set, $K$ is a skew field, $V$ is a left vector space over $K$, and $v_e\in V$ for each $e\in E$.  Then the set of linear dependencies among the vectors $v_e$, 
$D:=\{ X\in K^E: \sum_{e\in E} X_e v_e=0\}$,
is a left linear space over $K$. The collection of dependencies of minimal support
$$\mathcal{C}:=\{X\in D\del\{0\}: \text{if }Y\in D\del\{0\}\text{ and }\underline{Y}\subseteq \underline{X}, \text{ then }\underline{Y}=\underline{X}\},$$
satisfies the above left circuit axioms (C0)--(C3), so that $M(v_e: e\in E):=(E,\mathcal{C})$ is a left $K$-matroid. 


\subsection{The underlying matroid, circuit signatures, and coordinates}
If $M=(E,\mathcal{C})$ is a left- or right $H$-matroid, then $M$ determines an
{\em underlying} matroid $\underline{M}:=(E,\underline{\mathcal{C}})$, where
$$\underline{\mathcal{C}}:=\{\underline{X}: X\in \mathcal{C}\}.$$
If $H$ is the Krasner hyperfield, then $M$ in turn is uniquely determined by $\underline{M}$. Thus
a matroid $M$ over the Krasner hyperfield $\mathbb{K}$ is essentially a matroid. 

If $N$ is a matroid on $E$ and $H$ is a skew hyperfield, then a collection $\mathcal{C}\subseteq H^E$ is a {\em left $H$-signature} of $N$ if $\mathcal{C}$ satisfies (C0), (C1), and (C2), and 
$\underline{\mathcal{C}}$
is the collection of circuits of $N$. 

If $N$ is a matroid with bases $\mathcal{B}$, we name the set of ordered pairs of adjacent bases
$$A_N:=\{(B,B')\in \mathcal{B}\times\mathcal{B}: |B\setminus B'|=1\}.$$
Then a function $[.]: A_N\rightarrow H$ comprises {\em left $H$-coordinates} for $N$ if 
\begin{itemize}
\item [(CC0)] $[Fa, Fb]\x [Fb, Fa]=1$ if $Fa, Fb\in \mathcal{B}$.
\item [(CC1)] $[Fac,Fbc]\x [Fab, Fac]\x [Fbc, Fab]=-1$ if $Fab,  Fac, Fbc\in \mathcal{B}$.
\item [(CC2)] $[Fac, Fbc]=[Fad, Fbd]$ if $Fac, Fad, Fbc, Fbd\in \mathcal{B}$, but $Fab\not\in \mathcal{B}$.
\end{itemize}

As we will demonstrate, a left $H$-signature encodes the same information as left $H$-coordinates. If $\mathcal{C}$ is a left $H$-signature of $N$, then we may define a map $[.]: A_N\rightarrow H$ by setting 
$$[Fa, Fb]_{\mathcal{C}}:=-X_a^{-1}X_b$$
 where $X\in \mathcal{C}$ is any circuit such that $\underline{X}\subseteq Fab$. This is well-defined, since if $Y\in \mathcal{C}$ is any other circuit such that $\underline{Y}\subseteq Fab$, then $\underline{X}=\underline{Y}$ and hence by (C2) there exists an $\alpha\in H^\star$ so that $Y=\alpha X$. Then 
 $$Y_a^{-1}Y_b=(\alpha X_a)^{-1}(\alpha X_b)=X_a^{-1}X_b.$$
Conversely, given left coordinates $[.]$ for $N$, we put
$$\mathcal{C}_{N, [.]}:=\{X\in H^E: \underline{X}\text{ a circuit of }N
\text{ and }X_a^{-1}X_b=-[Fa, Fb]\text{ whenever }a,b\in \underline{X}\subseteq Fab\}.$$
We will usually omit the reference to $N$ when the choice of $N$ is unambiguous, and write $\mathcal{C}_{[.]}$.
\begin{lemma} \label{lem:coordinates}Let $N$ be a matroid on ground set $E$, let $\mathcal{C}\subseteq H^E$ and let $[.]:A_N\rightarrow H$. The following are equivalent.
\begin{enumerate}
\item $\mathcal{C}$ is a left $H$-signature of $N$, and $[.]=[.]_{\mathcal{C}}$.
\item $[.]$ are left $H$-coordinates, and $\mathcal{C}=\mathcal{C}_{[.]}$.
\end{enumerate}
\end{lemma}
\proof We show that (1) implies (2). Let $\mathcal{C}$ be a  left $H$-signature of $N$, and let $[.] = [.]_{\mathcal{C}}$. 
It suffices to show that the three axiom (CC0), (CC1), (CC2) hold for $[.]$.

(CC0): Note that if $Fa, Fb$ are both bases of $N$, and $X\in\mathcal C$ is any circuit so that $a,b\in \underline{X}\subseteq  Fab$, then 
$$[Fa,Fb][Fb,Fa]=(X_a^{-1}X_b)(X_b^{-1}X_a)=1.$$

(CC1): Assume $Fab, Fac, Fbc$ are bases of $N$. Then 
there exists a circuit $X\in \mathcal{C}$ so that $a,b,c\in \underline{X}\subseteq Fabc$. It follows that
$$[Fac,Fbc]\x [Fab, Fac]\x [Fbc, Fab]=-(X_a^{-1}X_b)(X_b^{-1}X_c)(X_c^{-1}X_a)=-1.$$

(CC2):  Assume that $Fac, Fad, Fbc, Fbd$ are bases of $N$. Then there are circuits $X, Y\in \mathcal{C}$, so that $a,b\in \underline{X}\subseteq Fabc$, and $a,b\in \underline{Y}\subseteq Fabd$.
If $Fab$ is not a basis of $N$, then $Fab$ contains a circuit, so that $\underline{X}=\underline{Y}$. By (C2), $Y=\alpha X$ for some $\alpha\in H^\star$. Then 
$$[Fac, Fbc]= -X_a^{-1}X_b=-(\alpha X_a)^{-1}(\alpha X_b)=-Y_a^{-1}Y_b=[Fad, Fbd].$$

We now argue that (2) implies (1). So suppose $[.]$ are left $H$-coordinates, and that $\mathcal{C}=\mathcal{C}_{[.]}$. We will first argue that for each circuit $C$ of $N$, there is an $X\in\mathcal{C}_{[.]}$ so that $\underline{X}=C$. 
So let $C$ be a circuit of $N$. 

Consider two elements $a,b\in C$. We claim that if $Fa, Fb, F'a, F'b$ are bases of $N$ so that if $C\subseteq Fab, F'ab$, then $[Fa,Fb]=[F'a, F'b]$.
To show this, we use induction on  $|F\setminus F'|$. Assume first that $|F\setminus F'|=1$. Then   $F=F''c$ and $F'=F''d$ for some $F'',c,d$. Since $C\subseteq (Fab)\cap(F'ab)=F''ab$, $F''ab$ is not a basis of $N$. By (CC2), it folllows that 
$$[Fa,Fb]=[F''ac, F''bc]=[F''ad, F''bd]=[F'a, F'b].$$
If $|F\setminus F'|>1$, then pick any $c\in (F\setminus F')\setminus C$. By virtue of the base exchange axiom in $N/C$, there exists a $d\in (F'\setminus F)\setminus C$ so that $F''a, F''b$ are  bases of $N$,  where $F''=F-c+d$. By the induction hypothesis, we obtain 
$$[Fa,Fb]=[F''a, F''b]=[F'a, F'b].$$
This proves the claim.

Fix any $c\in C$, let $B$ be a basis of $N$ containing $C-c$, and let $X\in H^E$ be such that $\underline{X}=C$, $X_c=1$, and $X_a:=[B-a+c, B]$ for all $a\in C-c$. By the claim, $X$ does not depend on the choice of $B$. By (CC0) and (CC1), we have 
$$X_a^{-1}X_b=(X_a^{-1}X_c)(X_c^{-1}X_b)=[Fab, Fbc][Fac, Fab]=-[Fac, Fbc]$$
whenever $a,b\in \underline{X}\subseteq Fabc$, so that  $X\in \mathcal{C}$. Thus 
$\underline{\mathcal{C}}$ is the set of circuits of $N$. It remains to  verify that $\mathcal{C}$ satisfies (C0), (C1), (C2), but these are straightforward.\endproof

The definition of right $H$-signatures $\mathcal{C}$, right coordinates $[.]$, and of the constructions $\mathcal{C}_{[.]}$ and $[.]_{\mathcal{C}}$ are obtained by reversing the order of multiplication throughout.

\subsection{The push-forward}
Let $f:H\rightarrow H'$ be a hyperfield homomorphism. Denote $f_* X:=(f(X_e): e \in E)$ for any $X\in H^E$, and for a set $\mathcal{C}\subseteq H^E$ denote
$$f_*\mathcal{C}:=\{\alpha'\x f_* X : \alpha'\in H', X\in \mathcal{C}\}.$$
From the definition of coordinates, it is immediate that $[B, B']_{f_* \mathcal{C}}=f([B,B']_{\mathcal{C}})$ for all adjacent bases $B,B'$.

If  $M=(E,\mathcal{C})$ is a left $H$-matroid, the {\em push-forward}  is $f_*M:=(E,f_*\mathcal{C})$. A straightforward verification yields that then (C0), (C1), (C2), (C3) hold for $f_*\mathcal{C}$, so that $f_* M$ is a left $H'$-matroid. 

Clearly $\underline{f_* M}=\underline{M}$ for any hyperfield homomorphism $f$ from $H$. In particular, if $\kappa: H\rightarrow \mathbb{K}$ then $\underline{M}=\underline{\kappa_* M}$, so that the underlying matroid can be considered as the ultimate push-forward. 

\subsection{Quasi-Pl\"ucker coordinates}
Let $H$ be a skew hyperfield and let $N$ be a matroid on $E$ with bases $\mathcal{B}$.  Then $[.]:A_N\rightarrow H$ are {\em left quasi-Pl\"ucker coordinates} if  
\begin{itemize}
\item[(P0)] $[Fa, Fb]\x [Fb, Fa]=1$ \hfill if $Fa, Fb\in \mathcal{B}$. 
\item[(P1)] $[Fac,Fbc]\x [Fab, Fac]\x [Fbc, Fab]=-1$ \hfill if $Fab,  Fac, Fbc\in \mathcal{B}$. 
\item[(P2)] $[Fa,Fb]\x [Fb, Fc]\x [Fc, Fa]=1$ \hfill if $Fa,  Fb, Fc\in \mathcal{B}$. 
\item[(P3)] $[Fac, Fbc]=[Fad, Fbd]$ \hfill if $Fac, Fad, Fbc, Fbd\in \mathcal{B}$, and $Fab\not\in \mathcal{B}$ or $Fcd\not\in \mathcal{B}$.
\item[(P4)] $1\in [Fbd, Fab]\x [Fac, Fcd]\+[Fad, Fab]\x [Fbc, Fcd]$ \hfill if $Fac, Fad, Fbc, Fbd, Fab, Fcd\in \mathcal{B}$.
\end{itemize} 
 We will show that in the presence of an underlying matroid $N$, these axioms are cryptomorphic to the left circuit axioms (C0)-(C3). 

\begin{theorem} \label{thm:plucker}Let $N$ be a matroid on $E$, let $H$ be a skew hyperfield, let $[.]:A_N\rightarrow H$ map, and let $\mathcal{C}\subseteq H^E$.
The following are equivalent:
\begin{enumerate}
\item $M=(E, \mathcal{C})$ is a  left $H$-matroid such that $\underline{M}=N$, and $[.]= [.]_{\mathcal{C}}$.
\item $[.]$ are left quasi-Pl\"ucker coordinates for $N$, and $\mathcal{C}=\mathcal{C}_{[.]}$.
\end{enumerate}
\end{theorem}
\proof We show that (1) implies (2). Let $M=(E,\mathcal{C})$ be a  left $H$-matroid such that $N=\underline{M}$, and let $[.] = [.]_{\mathcal{C}}$. By Lemma \ref{lem:coordinates}, $[.]$ are coordinates for $N$. We must show that the five axioms (P0)-(P4) hold. But (P0) is (CC0),  (P1) is (CC1), and (P3) partially follows from (CC2). We verify what remains.

(P2): Suppose that $Fa, Fb, Fc$ are bases of $N$, then
there are circuits $X,Y,Z\in \mathcal{C}$ so that $a,b\in \underline{X}\subseteq Fab$, and $b,c\in \underline{Y}\subseteq Fbc$, and $a,c\in \underline{Z}\subseteq Fac$ which determine the quasi-Pl\"ucker coordinates
$$[Fa, Fb]=-X_a^{-1}X_b, ~ [Fb, Fc]:=-Y_b^{-1}Y_c, ~[Fc, Fa]=-Z_c^{-1}Z_a.$$
The circuits $X,Y$ are modular, and by (C1) we may assume without loss of generality that $X_b=-Y_b$. By (C3) there exists a circuit $Z'\in \mathcal{C}$ with $Z'_b=0$ and $Z'\in X\+Y$. Then $\underline{Z'}\subseteq Fac$, so that $\underline{Z'}=\underline{Z}$. By (C2), we may assume that $Z=Z'$.
Then $Z_a=Z'_a\in X_a\+ Y_a=X_a\+0=\{X_a\}$ and $Z_c=Z'_c\in X_c\+ Y_c=0\+Y_c=\{Y_c\}$, so that 
$Z_a=X_a$ and $Z_c=Y_c$. It follows that
$$[Fa,Fb]\x [Fb, Fc]\x [Fc, Fa]=-(X_a^{-1}X_b)({Y_b}^{-1}Y_c)({Z_c}^{-1}Z_a)=1.$$

(P3): Assume that $Fac, Fad, Fbc, Fbd$ are bases of $N$. The case that $Fab\not\in\mathcal{B}$ is settled by (CC2), and we assume $Fcd\not\in \mathcal{B}$. Then there are circuits $X, Y\in \mathcal{C}$, so that $a,b\in \underline{X}\subseteq Fabc$, and $a,b\in \underline{Y}\subseteq Fabd$, and we may assume that $X_a=-Y_a$ by (C2). By (C3), there is a circuit $Z\in \mathcal{C}$ so that $\underline{Z}\subseteq Fbcd$, and $Z\in X\+ Y$. As $Fcd$ is dependent, we have $\underline{Z}\subseteq Fcd$, so that $0=Z_b\in X_b\+ Y_b$, i.e. $X_b=-Y_b$. Then also 
$$[Fac, Fbc]= -X_a^{-1}X_b=-Y_a^{-1}Y_b=[Fad, Fbd].$$

(P4): Assume that $Fac, Fad, Fbc, Fbd, Fab, Fcd$ are all bases of $N$. Then there are circuits $X, Y\in \mathcal{C}$, so that $a,c,d\in \underline{X}\subseteq Facd$, and $b,c,d\in \underline{Y}\subseteq Fbcd$.
Then $X,Y$ are modular, and by (C2) we may assume that $X_c=-Y_c$. By (C3), there is a circuit $Z\in\mathcal{C}$ so that $Z_c=0$ and $Z\in X\+ Y$. 
Thus $Z_a\in  X_a\+Y_a=X_a\+0=\{X_a\}$, $Z_b\in X_b\+ Y_b=0\+Y_b=\{Y_b\}$, and $Z_d\in X_d\+Y_d$. 
It follows that $Z_a=X_a\neq 0$, $Z_b=Y_b\neq 0$, and thus $a,b \in \underline{Z}\subseteq Fabd$. 
Since $Fab$  is a basis of $\underline{M}$, we have $\underline{Z}\not \subseteq Fab$, and hence $Z_d\neq 0$.  Then $$Z_d\in X_d\+Y_d=Z_aX_a^{-1}X_d \+Z_bY_b^{-1}Y_d$$
Multiplying on the left by $Z_d^{-1}$ and using the left distributivity of the hyperring $H$, it follows that
$$1 = Z_d^{-1}Z_d\in (Z_d^{-1}Z_a)\x(X_a^{-1}X_d) \+(Z_d^{-1}Z_b)\x(Y_b^{-1}Y_d)=[Fbd, Fab]\x [Fac, Fcd]\+[Fad, Fab]\x [Fbc, Fcd].$$
This completes the proof of (1)$\Rightarrow$(2).

We next show that (2) implies (1). Let $[.]:A_N\rightarrow H$ be left quasi-Pl\"ucker coordinates for $N$, and suppose that  $\mathcal{C}=\mathcal{C}_{[.]}$. By Lemma \ref{lem:coordinates}, $\mathcal{C}$ satisfies (C0), (C1), (C2).  It remains to show (C3).

So let $X, Y\in \mathcal{C}_{[.]}$ be modular, and consider a $c\in \underline{X}\cap \underline{Y}$. Assume that $X_c=-Y_c$. 
There exists a circuit $Z\in \mathcal{C}_{[.]}$ with $\underline{Z}\subseteq \underline{X}\cup\underline{Y}-c$, and we may assume that $Z_a=X_a$ for some $a\in \underline{X}\setminus \underline{Y}$. It remains to show that $Z\in X\+ Y$.

Pick $b\in \underline{Y}\setminus\underline{X}$. There is an $F$ such that $a,c\in \underline{X}\subseteq Fac$,  $b,c\in \underline{Y}\subseteq Fbc$, and $a,b\in \underline{Z}\subseteq Fab$. Then by (P2)
$$Z_a^{-1}Z_b=-[Fa,Fb]=-[Fa,Fc][Fc,Fb]=(X_a^{-1}X_c)(Y_c^{-1}Y_b)=X_a^{-1}Y_b,$$
so that $Z_b=Y_b\in 0\+Y_b=X_b\+Y_b$, as required.

Next, consider a $d\in \underline{X}\cap\underline{Y}$, other than $c$. We may assume that $Z_d=1$, again by rescaling as in (C1). By rescaling $X$ and $Y$ accordingly, we may assume that $X_a=Z_a$, $Y_b=Z_b$, and $X_c=-Y_c$. Then
$$X_d=-X_a[Fac, Fad]=-Z_a[Fac, Fcd]=Z_d[Fbd, Fab][Fac, Fad]=[Fbd, Fab][Fac, Fad]$$
and
$$Y_d=-Y_b[Fbc, Fcd]=-Z_b[Fbc, Fcd]=Z_d[Fad,Fab][Fbc, Fcd]=[Fad,Fab][Fbc, Fcd].$$
Hence by (P4),  $Z_d= 1\in [Fbd, Fab][Fac, Fad]\+[Fad,Fab][Fbc, Fcd]=X_d\+Y_d$. 
\endproof
\ignore{A {\em right Pl\"ucker map} is determined by the same axioms 
(P0)-(P4), but with the reverse order of factors in each product. A right $H$-matroid then determines a right Pl\"ucker map determined by   
$$[Fa, Fb]_M:=X_bX_a^{-1},$$
and conversely, a right Pl\"ucker map determines a right $H$-matroid with circuits
$$\mathcal{C}_{[.]}:=\{X\in H^E: \underline{X}\text{ is a circuit of }N,  \text{ and }X_bX_a^{-1}=[Fa, Fb]\text{ whenever }a,b\in \underline{X}\subseteq Fab\}.$$
}

\subsection{Duality} Let $H$ be a skew hyperfield, and let $E$ be a finite set. We say that $X, Y\in H^E$ are {\em orthogonal}, denoted $X\perp Y$, if 
$$0\in {\+}_{e\in E} X_e\x Y_e.$$
For sets $\mathcal{C},\mathcal{D}\subseteq H^E$, we write $\mathcal{C}\perp_k \mathcal{D}$ if $X\perp Y$ for all $X\in\mathcal{C}$ and $Y\in \mathcal{D}$ such that $|\underline{X}\cap\underline{Y}|\leq k$. 

Let $N$ be a matroid on $E$ and let $H$ be a skew hyperfield. To any $[.]:A_N\rightarrow H$ we associate a {\em dual map} 
$[.]^*:A_{N^*}\rightarrow H$ by setting
$$[B, B']^*:= - [E\setminus B, E\setminus B']$$
for all $(B, B')\in A_{N^*}$. It is evident from this definition that $[.]^{**}=[.]$. 
\begin{lemma} \label{lem:orth2} Let $N$ be a matroid on $E$ and let $H$ be a skew hyperfield, let $\mathcal{C}$ be a left $H$-signature of $N$, and let $\mathcal{D}\subseteq H^E$. The following are equivalent.
\begin{enumerate}
\item $\mathcal{D}$ is a right $H$-signature of $N^*$, and $\mathcal{C}\perp_2\mathcal{D}$
\item $[.]:=[.]_{\mathcal{C}}$ satisfies (P0), (P1), (P2), (P3), and $\mathcal{D}=\mathcal{C}_{[.]^*}$
\end{enumerate}
\end{lemma}
\proof We show that (1) implies (2). If $\mathcal{D}$ is a right $H$-signature of $N^*$, and $\mathcal{C}\perp_2\mathcal{D}$, then $[.]_{\mathcal{D}}=[.]_{\mathcal{C}}^*$. By Lemma \ref{lem:coordinates}, it follows that $\mathcal{D}=\mathcal{C}_{[.]^*}$.
Being right $H$-coordinates, $[.]_{\mathcal{C}}^*$  satisfies (CC0), (CC1), (CC2), which in terms of $[.]:=[.]_{\mathcal{C}}$ translates to 
\begin{itemize}
\item[(CC0)$^*$] $[Fa, Fb]\x [Fb, Fa]=1$  if $Fa, Fb\in \mathcal{B}$. 
\item[(CC1)$^*$] $[Fa,Fb]\x [Fb, Fc]\x [Fc, Fa]=1$ if $Fa,  Fb, Fc\in \mathcal{B}$. 
\item[(CC2)$^*$] $[Fac, Fbc]=[Fad, Fbd]$ if $Fac, Fad, Fbc, Fbd\in \mathcal{B}$, and  $Fcd\not\in \mathcal{B}$.
\end{itemize} 
Together with (CC0), (CC1), (CC2) for $[.]$, we have (P0), (P1), (P2), (P3) for $[.]$.

The proof that (2) implies (1) is a reversal of these steps.
\endproof
\begin{lemma}\label{lem:orth3}Let $N$ be a matroid on $E$ and let $H$ be a skew hyperfield, let $\mathcal{C}$ be a left $H$-signature of $N$, and let $\mathcal{D}\subseteq H^E$. The following are equivalent.
\begin{enumerate}
\item $\mathcal{D}$ is a right $H$-signature of $N^*$, and $\mathcal{C}\perp_3\mathcal{D}$.
\item $[.]:=[.]_{\mathcal{C}}$ are left quasi-Pl\"ucker coordinates, and $\mathcal{D}=\mathcal{C}_{[.]^*}$.
\end{enumerate}
\end{lemma}
\proof In view of Lemma \ref{lem:orth2}, we need to argue that if $\mathcal{C}$ is a left $H$-signature of $N$ and $\mathcal{D}$ is a right $H$-signature of $N^*$ so that $\mathcal{C}\perp_2\mathcal{D}$, then 
$$\mathcal{C}\perp_3\mathcal{D}\text{ if and only if (P4) holds for }[.]:=[.]_{\mathcal{C}}.$$
We first show sufficiency. So assume that $\mathcal{C}\perp_3\mathcal{D}$, and let $Fac, Fad, Fbc, Fbd, Fab, Fcd$ be bases of $N$. Let $X\in\mathcal{C}$ be such that $a,b,d\in \underline{X}\subseteq Fabd$, and let $Y\in\mathcal{D}$ be such that $a,b,d\in \underline{Y}\subseteq E\setminus Fd$. Without loss of generality, we may assume that $X_d=1$ and $Y_d=-1$. Using that $X\perp Y$, 
$$0\in X_a\cdot Y_a\+ X_b\cdot Y_b \+ X_d\cdot Y_d = [Fbd, Fab]\x [Fac, Fcd]\+[Fad, Fab]\x [Fbc, Fcd]\+ -1,$$
and it follows that $1\in [Fbd, Fab]\x [Fac, Fcd]\+[Fad, Fab]\x [Fbc, Fcd]$. 

To see necessity, let $X\in\mathcal{C}$ and $Y\in\mathcal{D}$ be such that $\underline{X}\cap \underline{Y}=\{a,b,d\}$ for distinct  $a,b,d\in E$. Since $\underline{Y}-ab$ is independent in $N^*$, we have $r(N\del (\underline{Y}-ab))=r(N)$. Hence, there exist a basis $Fab$ of $N$ extending the independent set $\underline{X}-d$ of $N$, such that $F\cap \underline{Y}=\emptyset$. By a dual argument, there exists a basis $Gab$ of $N^*\del F$ extending $\underline{Y}-d$. Since $|Fab|+|Gab|=r(N)+r(N^*)=|E|$ and $|F\cap G|=\emptyset$, $E\del(Fab\cup Gab)$ contains an element $c$ besides $d$. Scaling, we may assume that  $X_d=1$ and $Y_d=-1$. Using (P4), we have 
$$X_a\cdot Y_a\+ X_b\cdot Y_b \+ X_d\cdot Y_d = [Fbd, Fab]\x [Fac, Fcd]\+[Fad, Fab]\x [Fbc, Fcd]\+ -1\ni 0,$$
so that $X\perp Y$. 
\endproof
We say that a left  $H$-matroid $M=(E,\mathcal{C})$ and a right $H$-matroid $M'=(E, \mathcal{D})$ are {\em dual} if $\underline{M}=\underline{M'}^*$ and $\mathcal{C}\perp_3\mathcal{D}$. By Lemma \ref{lem:orth3}, each left or right $H$-matroid $M$ has a dual, which we denote by $M^*$. We highlight the following direct consequence of  Lemma \ref{lem:orth3}, using Theorem \ref{thm:plucker}.
\begin{theorem}\label{thm:duality} Let $N$ be a matroid on $E$ and let $H$ be a skew hyperfield. If $\mathcal{C}$ is a left $H$-signature of $N$ and $\mathcal{D}$ is a right $H$-signature of $N^*$ so that $\mathcal{C}\perp_3\mathcal{D}$, then $M=(E,\mathcal{C})$ is a left $H$-matroid and $M^*=(E,\mathcal{D})$.
\end{theorem}

\subsection{Minors} Let $N$ be a matroid on $E$, and let $\mathcal{C}$ be a left $H$-signature of $N$. For any disjoint sets $S,T\subseteq E$, 
put
$$\mathcal{C}/ S\del T:=\{X_{|E\del(S\cup T)}: X\in \mathcal{C}, \underline{X}\subseteq E\del T,\text{ and }\underline{X}\del S\text{ a circuit of }N/S\del T\}.$$
If $M=(E, \mathcal{C})$, the {\em minor} of $M$ obtained by {\em contracting} $S$ and {\em deleting } $T$ is $M/S\del T:=(E\del(S\cup T),  \mathcal{C}/ S\del T)$.
By construction, this minor  $M/S\del T$ is a left $H$-signature of $N/S\del T$.

For left coordinates $[.]:A_N\rightarrow H$ with associated circuit signature $\mathcal{C}:=\mathcal{C}_{[.]}$, we define $[.]/ S\del T:=[.]_{\mathcal{C}/ S\del T}$ for any pair of disjoint sets $S, T\subseteq E$. In the special case that $S$ is independent in $N$ and $T$ is independent in $N^*$, we have 
$$[Fa, Fb]/S\del T = -X_a^{-1}X_{b} = [S\cup Fa, S\cup Fa]$$
for any pair of adjacent bases $Fa, Fb$ of the minor $N/S\del T$ of $N$, where $X\in \mathcal{C}$ is any circuit so that $a,b\in \underline{X}\subseteq S\cup Fab$. 
Note that for any disjoint sets $S, T\subseteq E$, there exist disjoint sets $S', T'\subseteq E$ so that $N/S\del T=N/S'\del T'$, with $S'$ independent in $N$ and $T'$ is independent in $N^*$. 

\begin{lemma} \label{lem:minor} Let $N$ be a matroid on $E$, let $H$ be a skew hyperfield, and let $[.]:A_N\rightarrow H$ be $H$ coordinates for $N$. The following are equivalent.
\begin{enumerate}
\item $[.]$ are left quasi-Pl\"ucker coordinates.
\item $[.]/S\del T$ are left quasi-Pl\"ucker coordinates for all disjoint $S, T\subseteq E$ so that \begin{enumerate}
\item $S$ is independent in $N$ and $T$ is independent in $N^*$; and
\item $N/S\del T$ has rank $\leq 2$ and corank $\leq 2$.   
\end{enumerate}
\end{enumerate}
\end{lemma}
\proof That (1) implies (2) is straightforward. We prove that (2) implies (1). Assume (1). To see that (P0) holds for $[.]$, let $Fa, Fb$ be a bases. Then 
$$[Fa,Fb]\cdot [Fb, Fa]=[a, b]/S\del T \cdot [b, a]/S\del T=1$$
by (P0) for $[.]/S\del T$, where $S=F$ and $T=E\del Fab$. The minor $N/S\del T$ has ground set $E'=E\del (S\cup T)=ab$ and a basis $B'=Fa\del S=\{a\}$. Then the rank of $N/S\del T$ is $|B'|=1$ and the corank of $N$ is $|E'|-|B'|=1$.

An analogous argument applies to each of the other axioms: contract $S=F$ and delete $T=E\del Fabc$ (for (P1), (P2)) or $T=E\del Fabcd$ (for (P3), (P4)). In each case, the minor $N/S\del T$ has both rank $\leq 2$ and corank $\leq 2$.
\endproof
Using Theorem \ref{thm:plucker} to translate back to circuit signatures, we obtain: 
\begin{theorem}Let $N$ be a matroid on $E$, let $H$ be a skew hyperfield, and let $\mathcal{C}$ be a left $H$-signature of $N$. Then $M=(E, \mathcal{C})$ is a left $H$-matroid if and only if $M/S\del T$ is a left $H$-matroid, for all $S, T\subseteq E$ so that $N/S\del T$ has both rank and corank $\leq 2$. 
\end{theorem}
This theorem is known for valuated matroids, deriving from more general statements about {\em matroids over perfect fuzzy rings} due to Dress and Wenzel \cite[Section 3]{DressWenzel1992b}.

\subsection{The weak order} Let $M=(E, \mathcal{C}), M'=(E,\mathcal{C})$ be left $H$-matroids. We say that $M'$ is a {\em weak image} of $M$, notation $M'\preceq M$, if for all $X\in \mathcal{C}$ there exists an $X'\in \mathcal{C}'$ so that $X_e=X'_e$ for all $e\in \underline{X'}$. It follows that if $M'\preceq M$, then $f_* M'\preceq f_* M$ for any hyperfield homomorphism from $H$, and in particular that $\underline{M'}$ is a weak image of $\underline{M}$ in the usual sense for matroids.

If $\underline{M'}\preceq \underline{M}$ and $r(M')=r(M)$, then each basis of $M'$ is necessarily a basis of $M$. In this case,  we have $M'\preceq M$ if and only if $[B, B']_{M'}=[B, B']_{M}$ for all adjacent bases $B, B'$ of $M'$, i.e. if $[.]_{M'}$ is the restriction of $[.]_{M}$ to $A_{\underline{M'}}$.
\begin{lemma} \label{lem:weakimage}Let $M$ be a left $H$-matroid and let $N$ be a matroid, so that $N$ is a rank-preserving weak image of $\underline{M}$. Let $[.]:A_N\rightarrow H$ be the restriction of $[.]_M$ to $A_N$. Then $[.]$ are quasi-Pl\"ucker coordinates for $N$ if and only if $[.]$ satisfies (P3).\end{lemma}
\proof As $M$ is a left $H$-matroid,   (P0), (P1), (P2), and (P4) hold for $[.]_M$. The premise of each of these axioms is purely that certain bases exist. Since each basis of $N$ is necessarily a basis of $\underline{M}$, the same axioms will hold true for the restriction $[.]$ of $[.]_M$.  Hence if (P3) also holds for $[.]$, then $[.]$ are quasi-Pl\"ucker coordinates.\endproof

\subsection{Rescaling} If $N$ is a matroid on $E$, $\mathcal{C}$ is a left $H$-signature of $N$, and $\rho: E\rightarrow H^\star$, then {\em rescaling $\mathcal{C}$ by $\rho$} yields 
$$\mathcal{C}^\rho:=\{(X_e\rho_e: e\in E): X\in \mathcal{C}\}.$$
If $\mathcal{D}$ is a right $H$-signature of $N^*$, and $\rho: E\rightarrow H^\star$, then {\em rescaling $\mathcal{D}$ by $\rho$}
yields 
$$\mathcal{D}^\rho:=\{(\rho_e Y_e: e\in E): Y\in \mathcal{D}\}.$$
\begin{lemma} \label{lem:rescaling}Let  $N$ be a matroid on $E$, let $\mathcal{C}$ be a left $H$-signature of $N$ and let $\mathcal{D}$ is a right $H$-signature of $N^*$. Then 
for any $\rho:E\rightarrow H^\star$ we have $\mathcal{C}\perp_k \mathcal{D}\text{ if and only if }\mathcal{C}^{\rho^{-1}}\perp_k \mathcal{D}^{\rho}$,
where $\rho^{-1}: e\mapsto \rho_e^{-1}$.
\end{lemma}
For a left or right $H$-matroid $M$ on $E$ with circuits $\mathcal{C}$ and cocircuits $\mathcal{D}$, {\em rescaling $M$ by $\rho:E\rightarrow H^\star$} yields a matroid $M^\rho$ with circuits $\mathcal{C}^{\rho^{-1}}$ and cocircuits $\mathcal{D}^{\rho}$. 
This rescaling convention (as opposed to scaling $\mathcal{C}$ by $\rho$ and $\mathcal{D}$ by $\rho^{-1}$) is consistent with the effect of scaling vectors $v_e$ from a left vector space $V$: we have $M(\rho_ev_e: e\in E)=M(v_e: e\in E)^\rho$.
If $M$ is a left $H$-matroid, then for the coordinates of $M^\rho$ we obtain 
$$[Fa, Fb]_{M^\rho}=   -(X_a\rho_a^{-1})^{-1}(X_b\rho_b^{-1}) =-\rho_aX_a^{-1}X_b\rho_b^{-1} = \rho_a[Fa, Fb]_{M}\rho_b^{-1}$$
for any $X\in\mathcal{C}$ so that $a,b\in \underline{X}\subseteq Fab$. 
For a right $H$-matroid $M$, we have a reversed order of multiplication: $[Fa, Fb]_{M^\rho}=\rho_b^{-1}[Fa, Fb]_{M}\rho_a$.

\ignore{Rescaled signatures $\mathcal{C}^\rho$ and $\mathcal{D}^{\rho}$ have coordinates 
$$[Fa, Fb]_{\mathcal{C}^\rho} =   (X_a\rho_a)^{-1}(X_b\rho_b) =\rho_a^{-1}X_a^{-1}X_b\rho_b = \rho_a^{-1}[Fa, Fb]_{\mathcal{C}}\rho_b$$
for any $X\in\mathcal{C}$ so that $a,b\in \underline{X}\subseteq Fab$, and
$$[Fa, Fb]_{\mathcal{D}^\rho} =   (\rho_bY_b)(\rho_aY_a)^{-1} =\rho_bY_bY_a^{-1} \rho_a^{-1}= \rho_b[Fa, Fb]_{\mathcal{D}}\rho_a^{-1}$$
for any $Y\in\mathcal{D}$ so that $a,b\in \underline{Y}\subseteq Fab$.
We define {\em rescaling} of left and right coordinates by $\rho$ accordingly, setting $[Fa, Fb]^\rho:=\rho_a^{-1}[Fa, Fb]\rho_b$ for left coordinates and $[Fa, Fb]^\rho:=\rho_b[Fa, Fb]\rho_a^{-1}$ for right coordinates. 
}
We say that $\mathcal{C}$ and $\mathcal{C}'$ are {\em rescaling equivalent} if $\mathcal{C}^\rho=\mathcal{C}'$ for some $\rho:E\mapsto H^\star$, and write $\mathcal{C}\sim \mathcal{C}'$ .
We investigate the rescaling classes of $U_{2,4}$. For any $x,y\in H^\star$, let   $\mathcal{U}_H(x,y)$ denote the unique $H$-signature of $U_{2,4}$ containing  $$(0,1,1,1), (1,0,-1,-x), (1,1,0,y), (1,x,-y,0).$$
 \begin{lemma} Let $H$ be a skew hyperfield, and let $M=(E,\mathcal{C})$ be a left $H$-matroid so that  $\underline{M}=U_{2,4}$. 
 Then there are  $x,y\in H^\star$ with $1\in x\+ y$ so that  $\mathcal{C}\sim\mathcal{U}_H(x,y)$. Moreover, $$\{(x',y'):  \mathcal{C}\sim\mathcal{U}_H(x',y')\}=\{\beta(x, y)\beta^{-1}:\beta\in H^\star\}.$$
\end{lemma} 
\proof Write $E=\{a,b,c,d\}$, and pick $W, X, Y, Z\in \mathcal{C}$ such that $\underline{W}=bcd, \underline{X}=acd, \underline{Y}=abd, \underline{Z}=abc$. 
Using (C2), we may assume that $X_a=Y_a=Z_a=1$, and $W_b=Y_b$. Define $\rho\in E\rightarrow H^\star$ by 
$$\rho_a=1, \rho_b=W_b^{-1}, \rho_c=W_c^{-1},   \rho_d=W_d^{-1}.$$ 
Replacing $\mathcal{C}$ with $\mathcal{C}^\rho\sim \mathcal{C}$, we have $W=(0,1,1,1), X=(1,0,s,-x), Y=(1,1,0,y), Z=(1,x',-y',0)$ for some $s,x,y,x',y'\in H^\star$. Note that each pair of these circuits is modular in $\mathcal{C}$. Applying (C3), we have
\begin{enumerate}
\item $X\in (-W)\+ Y$, so that $s=X_c\in (-W_c)\+Y_c=\{-1\}$, so $s=-1$;  
\item $Z\in xW\+ X$, so that $x'= Z_b\in xW_b\+X_b=x\+ 0$, so  $x'=x$;
\item $Z\in (-yW)\+ Y$, so that $-y'= Z_c\in (-yW_c)\+Y_c=-y\+ 0$, so  $y'=y$; and 
\item $W\in (-X)\+Y$, so that $1=W_d\in (-X_d)\+ Y_d=x\+y$.
\end{enumerate}
Then $\mathcal{C}=\mathcal{U}_H(x,y)$ and $1\in x\+ y$, as required. Finally, if $\mathcal{U}_H(x',y')\sim\mathcal{U}_H(x,y)$, then we must have $\mathcal{U}_H(x',y')=\mathcal{U}_H(x,y)^\rho$ with $\rho=\beta\one_E$ for some $\beta$. It then follows that $(x', y')=(\beta x\beta^{-1}, \beta y\beta^{-1})$.
\endproof
Thus the conjugacy class of the pair $(x,y)$ as in the lemma is a scaling invariant of any $H$-orientation of $U_{2,4}$, and more generally, gives an invariant for each $U_{2,4}$-minor of each left $H$-matroid $M$. 
\subsection{\label{ss:cr} Cross ratios}
Let $M$ be a left $H$-matroid on $E$. The {\em cross ratio} is defined as 
 $$cr_M(F,a,b,c,d):=[Fac,Fad]_M[Fbd, Fbc]_M.$$
Formally $cr_M: CR_{\underline{M}}\rightarrow H$, where $CR_N:=\{(F,a,b,c,d): Fac, Fad, Fbd, Fbc\text{ are bases of }N\}.$ The following properties follow by substituting the definition of cross ratio and applying the quasi-Pl\"ucker axioms.
\begin{itemize}
\item[(CR0)] $cr(F,a,b,c,d)cr(F,b,a,c,d)=1$.
\item[(CR1)] $cr(F,a,b,d,e)cr(F,b,c,d,e)cr(F,c,a,d,e)=1$.
\item[(CR2)] $cr(Fa,b,c,d,e)cr(Fc,a,b,d,e)cr(Fb,c,a,d,e)=1$.
\item[(CR3)] $cr(F,a,b,c,d)=1$ if $Fab$ or $Fcd$ is not a basis of $\underline{M}$.
\item[(CR4)] $1\in cr(F,b,c,d,a)\+cr(F,a,c,d,b)$.
\item[(CRP)] $[Fad,Fcd] cr(F,a,b,c,d)= cr(F,c,b,a,d) [Fbd,Fcd]$.
\end{itemize}
In the context of quasi-determinants, the cross ratio was similarly defined by Gelfand, Gelfand, Retakh, and Wilson, who also note such properties  \cite{GGRW2005, Retakh2014}.

\subsection{Matroids over commutative hyperfields} If $N$ is a matroid on $E$ of rank $r$ and $H$ is a hyperfield, then a {\em Grassmann-Pl\"ucker} function for $N$ is a function $\phi:E^r\rightarrow H$ such that
\begin{itemize}
\item[(GP0)] $\phi(B)\neq 0$ if and only if $\underline{B}$ is a basis of $N$.
\item[(GP1)] $\phi(B^\tau)=\mbox{sign}(\tau)\phi(B)$ for all $B\in E^r$ and permutations $\tau$ of $[r]$.
\item[(GP2)] $0\in \phi(Fab)\phi(Fcd)~\+~ \phi(Fac)\phi(Fdb)~\+ ~\phi(Fad)\phi(Fbc)$
for all $F\in E^{r-2}$ and $a,b,c,d\in E$.
\end{itemize}
In the above axioms and in the remainder of this section, we use the following notation. If $F\in E^k$, we denote the underlying set as $\underline{F}:=\{F_1,\ldots, F_k\}$,  $F^\tau:=(F_{\tau(1)},\ldots, F_{\tau(k)})$ for any permutation $\tau$ of $[k]:=\{1,\ldots,k\}$, and for any $a\in E$ we put $Fa:=(F_1,\ldots, F_k, a)\in E^{k+1}$ .

Grassmann-Pl\"ucker functions are closely related to quasi-Pl\"ucker coordinates. The proof of the following lemma amounts to a straightforward verification, which we omit.
\begin{lemma} \label{lem:gp21}Let $N$ be a matroid and let $H$ be a commutative hyperfield. Suppose $\phi:E^r\rightarrow H$ is a Grassmann-Pl\"ucker function for $N$. There is a unique function $[.]:A_N\rightarrow H$ such that $$[\underline{F}a, \underline{F}b]=\phi(Fa)/\phi(Fb)$$ for all $F\in E^{r-1}$ and $a,b\in E$ so that $\underline{F}a, \underline{F}b$ are bases of $N$. Such  $[.]$ are quasi-Pl\"ucker coordinates for $N$.
\end{lemma}

\newcommand{\lb}{\ldbrack}
\newcommand{\rb}{\rdbrack}

The {\em Tutte group} of a matroid $N$ with bases $\mathcal{B}$ was defined by Dress and Wenzel in \cite{DressWenzel1989} as the abelian group $\mathbb{T}_N$ with a generator $\epsilon$ subject to the relation $\epsilon^2=1$, and a generator $\lb B, B'\rb$ for each $(B, B')\in A_N$ satisfying further relations 
\begin{itemize}
\item[(T0)] $\lb Fa, Fb\rb  \lb Fb, Fa\rb =1$ \hfill if $Fa, Fb\in \mathcal{B}$. 
\item[(T1)] $\lb Fac,Fbc\rb  \lb Fab, Fac\rb \lb Fbc, Fab\rb =\epsilon$ \hfill if $Fab,  Fac, Fbc\in \mathcal{B}$. 
\item[(T2)] $\lb Fa,Fb\rb  \lb Fb, Fc\rb  \lb Fc, Fa\rb =1$ \hfill if $Fa,  Fb, Fc\in \mathcal{B}$. 
\item[(T3)] $\lb Fac, Fbc\rb =\lb Fad, Fbd\rb $ \hfill if $Fac, Fad, Fbc, Fbd\in \mathcal{B}$, and $Fab\not\in \mathcal{B}$ or $Fcd\not\in \mathcal{B}$.
\end{itemize}
A comparison with the quasi-Pl\"ucker axioms (P0) --- (P3) immediately gives the following.
\begin{lemma} \label{lem:th}Let $N$ be a matroid, let $H$ be a commutative hyperfield, and let $[.]: A_N\rightarrow H$ be a function satisfying (P0), (P1), (P2), and (P3). There is a group homomorphism $h: \mathbb{T}_N\rightarrow H^\star$ so that $h: \epsilon \mapsto -1$ and $h: \lb B, B'\rb \mapsto [B, B']$ for all $(B, B')\in A_N$.\end{lemma}
In \cite{DressWenzel1989}, Dress and Wenzel define several further abelian groups from a matroid $N$, and show that each group is  (essentially) isomorphic to $\mathbb{T}_N$.  There seems to be a close relation between each presentation of the Tutte group and different axiomatizations of matroids over commutative hyperfields $H$, which could be characterized as multiplicative group homomorphisms from the Tutte group to $H^\star$ satisfying a further additive duality constraint. With this in mind, we will use one of their isomorphisms here to argue the converse of Lemma \ref{lem:gp21}. 

Let  $\mathbb{T}_N^{\mathcal{B}}$ be the abelian group with a generator $\epsilon$ so that $\epsilon^2=1$, and generator $\lb B\rb$ for each $B\in E^r$ such that $\underline{B}$ is a basis of $N$, satisfying the relations
\begin{itemize}
\item[(TB1)] $\lb B^\tau\rb=\epsilon \lb B\rb$ \hfill whenever $\mbox{sign}(\tau)=-1$
 \item[(TB2)] $\lb Fac\rb \lb Fbc\rb^{-1}=\lb Fad\rb \lb Fbd\rb^{-1} $ \hfill if $Fac, Fad, Fbc, Fbd\in \mathcal{B}$, and $Fab\not\in \mathcal{B}$ or $Fcd\not\in \mathcal{B}$.
\end{itemize}
The following is a direct consequence of \cite[Theorem 1.1]{DressWenzel1989}.
\begin{lemma} \label{lem:tt}Let $N$ be a matroid on $E$ of rank $r$. There is a group homomorphism $t: \mathbb{T}_N^{\mathcal{B}}\rightarrow \mathbb{T}_N$ so that  $$t(\lb Fa\rb \lb Fb\rb^{-1})=\lb Fa, Fb\rb$$ for all $F\in E^{r-1}$ and $a,b\in E$ so that $\underline{F}a, \underline{F}b$ are both bases of $N$, and $t(\epsilon)=\epsilon$. \end{lemma}

\ignore{
The converse direction could be argued using a cryptomorphism between two presentations of the {\em Tutte group}, show by Dress and Wenzel \cite[Theorem 1.1]{DressWenzel1989}. Setting up both presentations of the Tutte group and interpreting their theorem in our context takes almost as much space as giving a complete proof. We prefer to do the latter even if our argument is essentially the same, using Maurer's Homotopy Theorem \cite{Maurer1973}.   
\begin{lemma} \label{lem:gp12} Let $N$ be a matroid on $E$ of rank $r$, let $H$ be a commutative hyperfield, and let $[.]: A_N\rightarrow H$ be a function satisfying (P0), (P1), (P2), and (P3). There exists a $\phi:E^r\rightarrow H$ satisfying (GP1) and (GP2), so that $\phi(Fa)/\phi(Fb) =[\underline{F}a, \underline{F}b]$ for all 
for all $F\in E^{r-1}$ and $a,b\in E$ so that $\underline{F}a, \underline{F}b$ are both bases of $N$. \end{lemma}
\proof Let $\mathcal{B}$ denote the set of bases of $N$. Consider the weighted directed graph $D$ with vertices $V(D)=\{B\in E^r:\underline{B}\in\mathcal{B}\}$ and arcs 
\begin{enumerate}
\item $(B, B^\tau)$ of weight $w(B, B^\tau):=\text{sign}(\tau)$, for each $B\in V(D)$ and each permutation $\tau$ of $[r]$; and
\item $(B, B')$ of weight $w(B, B'):=[\underline{B}, \underline{B'}]$ for $B, B'\in V(D)$ so that $B_i\neq B'_i$ for exactly one $i\in [r]$.
\end{enumerate}
The weight of a walk $(B_1,\ldots, B_k)$ is defined as $w(B_1,\ldots, B_k):=w(B_1, B_2)w(B_2, B_3)\cdots w(B_{k-1}, B_k)$. If we assume that each directed closed walk in $D$ has total multiplicative weight 1, we may construct a $\phi$ satisfying the conditions of the theorem as follows. Fix any $B_0\in V(D)$, and set the value of $\phi(B)$ to the weight of any directed walk from $B$ to $B_0$ in $D$ if $\underline{B}\in\mathcal{B}$, and to $0$ otherwise. By our assumption that each closed walk in $D$ has weight 1, the value of $\phi(B)$ does not depend on the directed path from $B$ to $B_0$ we choose to fix $\phi(B)$, and hence we have  $\phi(B)/\phi(B')=w(B, B')$ for each arc $(B, B')\in A(D)$. Then $\phi$ is alternating by the presence of arcs (1), and by the presence of arcs (2) we have $\phi(Fa)/\phi(Fb)=w(Fa, Fb)=[\underline{F}a, \underline{F}b]$, as required.

It remains to justify the assumption that each directed closed walk in $D$ has multiplicative weight 1. We will first establish this for each of the closed walks of the following restricted  types.
\begin{enumerate}
\item $(B, B)$ for each $B\in E^r$;
\item $(B,B^\tau,B)$ for each $B\in E^{r}$ and permutation $\tau$;
\item $(B,B^\tau,B'^\tau,B',B)$, for each $B, B'\in V(D)$ so that $B_i\neq B'_i$ for exactly one $i$, and each $\tau$;
\item $(Fa,Fb,Fc,Fa)$, for each $F\in E^{r-1}$ and $a,b,c\in E$ so that $\underline{F}a, \underline{F}b, \underline{F}c\in\mathcal{B}$;
\item $(Fab,Fcb,Fca,Fba,Fab)$, for each $F\in E^{r-2}$ and $a,b,c\in E$ so that $\underline{F}ab, \underline{F}bc, \underline{F}ac\in\mathcal{B}$;
\item $(Fac,Fad,Fbd,Fbc,Fac)$, for each $F\in E^{r-2}$ and $a,b,c,d\in E$ so that $\underline{F}ac, \underline{F}ad$, $\underline{F}bc, \underline{F}bd\in\mathcal{B}$, but $\underline{F}ab,\underline{F}cd\not\in\mathcal{B}$.
\end{enumerate}
For walks of type (1), we clearly have $w(B,B)=\mbox{sign}(id)=1$. A walk of type (2) has total weight
$w(B, B^\tau, B)=w(B, B^\tau)w(B^\tau, B)=\text{sign}(\tau)\text{sign}(\tau^{-1})=1$
and a similar argument applies to type (3). For type (4), we have 
$$w(Fa, Fb, Fc, Fa)=w(Fa, Fb)w(Fb, Fc)w(Fc, Fa)=[\underline{F}a, \underline{F}b][\underline{F}b, \underline{F}c][\underline{F}c, \underline{F}a]=1$$ by axiom (P2). Type (5) similary follows axiom (P1). For a  walk of type (6), we have
$$w(Fac,Fad,Fbd,Fbc,Fac)=[\underline{F}ac, \underline{F}ad][\underline{F}ad, \underline{F}bd]
[\underline{F}bd, \underline{F}bc][\underline{F}bc, \underline{F}ac].$$
Using the commutativity of $H$, this equals
$[\underline{F}ac, \underline{F}ad][\underline{F}bd, \underline{F}bc]\cdot
[\underline{F}ad, \underline{F}bd][\underline{F}bc, \underline{F}ac]=1$
by two applications of (P3).

Upon contracting all arcs of type (1) in $D$ and then taking the underlying simple graph, we obtain the base exchange graph $G$ of $N$, with vertices $V(G)=\mathcal{B}$ and edges $E(G)=\{\{B, B'\}: (B,B')\in A_{N}\}$. Each closed directed walk $W$ of $D$ induces an underlying closed walk $\underline{W}$ in $G$, in which each node $B$ in $W$ is matched by its underlying set $\underline{B}$.
Maurer's Homotopy Theorem \cite{Maurer1973} implies that each closed walk of $G$ can be reduced to a trivial closed walk by inserting and removing detours of the form $(Fa, Fb, Fa), (Fa, Fb, Fc, Fa), (Fab, Fbc, Fac, Fab)$, and $(Fac, Fad, Fbd, Fbc, Fac)$. If $Fab$ or $Fcd$ is a basis of $N$, then the quadrangle $(Fac, Fad, Fbd, Fbc, Fac)$ is not necessary for such reduction: the same reduction is obtained by using four triangles through $Fab$ or four triangles through $Fcd$. We may assume that each of the digons, triangles and quadrangles used for reduction in $G$ underly some walk of type (3), (4), (5), or (6) in $D$. 

For a directed closed walk $W$ in $D$, let $f_W:A(D)\rightarrow \Z$ denote the flow in $D$ induced by $W$, where $f_W(B,B')$ counts the number of time $W$ traverses the arc $(B, B')\in A(D)$. The multiplicative weight of a flow $f$ is $w(f):=\prod_{a\in A(D)} w(a)^{f(a)}$. Clearly $w(F_W)=w(W)$ for any closed walk $W$. 

Let $W$ be a directed closed walk in $D$, and let $W_1,\ldots, W_k$ be a collection of closed walks of types (3), (4), (5), and (6) so that inserting $\underline{W_1},\ldots,  \underline{W_l}$ and removing $\underline{W_{l+1}},\ldots, \underline{W_{k}}$ reduces $\underline{W}$ to a trivial walk in $G$. Defining $f:= f_W+\sum_{i=1}^l f_{W_i}-\sum_{i=l+1}^{k} f_{W_i}$, we have
$$w(f)=w(f_W) \left(\prod_{i=1}^lw(f_{W_i})\right)/\left(\prod_{i=l+1}^{k}w(f_{W_i})\right)=w(W),$$
since $w(f_{W_i})=w(W_i)=1$ for each $i$.
By construction of $f$, the total flow though each edge of $G$ is reduced to zero by adding and substracting the flows $f_{W_i}$, that is, for each $(A, A')\in A_M$, we have $$\sum \{ f(B, B')-f(B',B): \underline{B}=A, \underline{B'}=A'\} = 0.$$
It follows that there are walks $W_{k+1},\ldots, W_m$ of type (3) in $D$ so that $g:=f-\sum_{i=k+1}^mf(W_i)$ has $g(B,B')=0$ whenever $\underline{B}\neq \underline{B'}$. Clearly such a flow $g$  the sum of flows induced by walks of type (1) and (2) in $D$, and hence $w(W)=w(f)=w(g)=1$, as required. 
\endproof 
}
\begin{theorem}[Baker and Bowler \cite{BakerBowler2017}]\label{thm:gp} Let $H$ be a commutative hyperfield and let $M$ be an $H$-matroid on $E$. 
Then there exists a Grassmann-Pl\"ucker function $\phi$ for $\underline{M}$, so that 
$\phi(Fa)/\phi(Fb) =[\underline{F}a, \underline{F}b]_M$ 
for all $F\in E^{r-1}$ and $a,b\in E$ so that $\underline{F}a, \underline{F}b$ are both bases of $N$.
\end{theorem}
\proof 
Let $h$ and $t$ be group homomorphisms as in Lemma \ref{lem:th} and \ref{lem:tt}, repectively.
Let $\phi:E^r\rightarrow H$ be defined by setting $\phi(B)=0$ if $\underline{B}$ is not a basis of $N$ and $\phi(B)=h(f(\lb B\rb))$ otherwise. Then $\phi$ satisfies (GP0) and (GP1) by construction, and $\phi(Fa)/\phi(Fb) =[\underline{F}a, \underline{F}b]_M$ for all $F\in E^{r-1}$ such that $(\underline{F}a, \underline{F}b)\in A_{\underline{M}}$. It remains to show that $\phi$ satisfies (GP2). 

So let $F\in E^{r-2}$ and $a,b,c,d\in E$. If each of $\phi(Fab)\phi(Fcd), \phi(Fac)\phi(Fdb), \phi(Fad)\phi(Fbc)$ is nonzero, then by (P4) we have
$$1\in [\underline{F}bd, \underline{F}ab]_M[\underline{F}ac, \underline{F}cd]_M\+[\underline{F}ad, \underline{F}ab]_M [\underline{F}bc, \underline{F}cd]_M = 
\frac{\phi(Fbd)}{\phi(Fab)}\frac{\phi(Fac)}{\phi(Fcd)}\+\frac{\phi(Fad)}{\phi(Fab)}\frac{\phi(Fbc)}{\phi(Fcd)}$$
Multiplying by $\phi(Fab)\phi(Fcd)$ throughout, we obtain
$$\phi(Fab)\phi(Fcd)\in -\phi(Fac)\phi(Fdb)~\+ ~-\phi(Fad)\phi(Fbc)$$
which implies $0\in \phi(Fab)\phi(Fcd)~\+~\phi(Fac)\phi(Fdb)~\+ ~\phi(Fad)\phi(Fbc)$, as required. If one of the three terms is zero, we may assume by symmetry of 
$b,c,d$ that $\phi(Fab)\phi(Fcd)=0$, so that $\underline{F}ab$ or $\underline{F}cd$ is not a basis of $\underline{M}$. If one of $\phi(Fac)\phi(Fdb), \phi(Fad)\phi(Fbc)$ is nonzero, then so is the other, and then by (P3) we have 
$$\frac{\phi(Fac)}{\phi(Fbc)}=[\underline{F}ac, \underline{F}bc]_M= [\underline{F}ad, \underline{F}bd]_M=\frac{\phi(Fad)}{\phi(Fbd)}$$
which implies $0\in \phi(Fab)\phi(Fcd)~\+~\phi(Fac)\phi(Fdb)~\+ ~\phi(Fad)\phi(Fbc)$ upon cross multiplying.  \endproof
Over skew hyperfields, there seems to be no proper analogue of Grassmann-Pl\"ucker functions. However, with each skew hyperfield $H$ we may associate a commutative hyperfield $H^{ab}$, which arises by dividing out the commutator subgroup of $H^\star$, and there is a canonical homomorphism $\delta: H\rightarrow H^{ab}$. If $M$ is a left- or right $H$-matroid, then it may not be possible to define a Grassmann-Pl\"ucker function for $M$, but the push-forward $\delta_* M$ is a matroid over a commutative hyperfield, which does admit a Grassmann-Pl\"ucker function. 

\subsection{Quasi-determinants of matrices and matroids over skew fields}
For an invertible square $I\times J$ matrix $A$ with entries in a skew field $K$, Gelfand, Gelfand, Retakh, and Wilson  defined the {\em quasi-determinant } $|A|_{ij}:=b_{ji}^{-1}$, where $b_{ji}$ is the $ji$-th entry of the $J\times I$ matrix inverse $B$ of $A$ \cite[Definition 1.2.2]{GGRW2005}. There is a direct relation beween these quasi-determinants and quasi-Pl\"ucker coordinates of a right $K$-matroid arising from $A$. 

For an $r\times E$ matrix $V$ over $K$, we denote $M(V):=M(v_e: e\in E)$, where $v_e$ is the $e$-th column of $V$. We will consider the vectors $v_e$ as coming from a {\em right} vector space over $K$, so that $M(V)$ is defined as a right $K$-matroid. For any $F\subseteq E$, let $V[F]$ denote the restriction of $V$ to the columns indexed by $F$. 
\begin{lemma} \label{lem:qd}Let $A$ be an invertible $I\times J$ matrix over $K$, and let $V:=[I | A]$ be the $I\times (I\cup J)$ matrix so that $V[I]$ is the identity matrix and $V[J]=A$. 
 Then $[J, J-j+i]_{M(V)}=|A|_{ij}$.
\end{lemma}
\proof Let $B$ be the matrix inverse of $A$. Then $M(BV)=M(V)$, and we have $BV=B[I | A]=[B | I]$. For the circuit $X$ of $M(BV)=M(V)$ with $i,j\in \underline{X}\subseteq Ji$, and $X_i=1$, we have $X_j=-b_{ji}$, since 
$$b_{ji}X_i+X_j=(BV)_{ji}X_i+(BV)_{jj}X_j=\sum_e  (BV)_{je}X_e=0.$$
Then $[J, J-j+i]_{M(V)}=-X_iX_j^{-1}=b_{ji}^{-1}=|A|_{ij}$, as required.
\endproof
For a general $S\times T$ matrix $A$ with entries in $K$, the {\em left quasi-Pl\"ucker coordinate} is defined in  \cite{GGRW2005} as $q^I_{ij}(A):=|A[Ii]|_{si}^{-1}|\cdot A[Ij]|_{sj}$, after showing that this expression does not depend on the choice of $s\in S$. In terms of the matrix $V:= [I | A]$ and the right $K$-matroid $M=M(V)$, we have 
$$q^I_{ij}(A):=|A[Ii]|_{si}^{-1}\cdot|A[Ij]|_{sj}=[Ii, Is]^{-1}_M\cdot [Ij, Is]_M=[Ij, Ii]_M,$$ using Lemma \ref{lem:qd} and the multiplicative relation (P2) for right coordinates. 

Among the results in \cite[Section 4.4]{GGRW2005}, there are statements about the quasi-Pl\"ucker coordinates of a matrix $A$ corresponding to each of the axioms (P0)---(P4) we have used to define quasi-Pl\"ucker coordinates for matroids. In \cite{LR2017}, Laugwitz and Retakh consider an algebra $\mathcal{Q}_n$ with abstract generators $q_{ij}^I$ and defining relations similar to our axioms.

The  {\em Dieudonne determinant} \cite{Dieudonne1943} of an invertible  $I\times J$ matrix $A$ over $K$ equals $\phi(J)/\phi_M(I)$, where $\phi$ is any Grassmann-Plucker function for the matroid $M:=\delta_* M([I | A])$ and $\delta: K\rightarrow K^{ab}$ is the canonical hyperfield homomorphism to the abelianization of $K$.


\section{\label{sec:hyperfield}A skew hyperfield}
\subsection{\label{sigma}The skew hyperfield of monomials} Let $H$ be any hyperring, and let $\sigma:H\rightarrow H$ be an automorphism. 
We define a new hyperring $$H(T,\sigma,\min)=(\{T^\infty\}\cup \{a T^i: a\in H^\star, i\in\Z\}, 1, 0, \cdot, \boxplus),$$
as follows. As the notation suggests, we identify $a\in H$ with $aT^0$ and write $T^i$ for $1T^i$. We put $1:=T^0$ and $0:=T^\infty$. Multiplication follows the rules $0\cdot aT^i=aT^i\cdot 0=0$ and 
$$a T^i \cdot b T^j:=a\sigma^i(b)T^{i+j}$$
for all $a,b\in H^\star$ and $i,j\in\Z$. In particular, $a\cdot T^j=aT^j$. The hypersum is given by $0\+ x= x\+0=\{x\}$ and 
$$aT^i\boxplus bT^j:=\left\{\begin{array}{ll} 
\{aT^i\}&\text{if } i<j\\
\{bT^j\}&\text{if } i>j\\
(a+b)\cdot T^i&\text{if } i=j\text{ and } a\neq -b\\
(a+b)\cdot T^i \cup H^\star\cdot\{T^k: k\in \Z, k> i\}& \text{if }i=j\text{ and }a=-b
\end{array}
\right.
$$
for $a,b\in H^\star$ and $i,j\in\Z$, where $+$ is the hyperaddition of $H$. Note that in the last line of this definition, we have $0=0\cdot T^i\in (a+b)\cdot T^i$ as $a=-b$. 

There is a variant $H(T, \sigma, \max)$ which arises by reversing $<$ and $>$ in the above definition. In the present paper, we will hardly use this variant, and we will not substitute the symbol $T$. For brevity, we write $H^\sigma:=H(T,\sigma, \min)$ in what follows.
\begin{lemma}Let $H$ be a hyperring, and let $\sigma$ be an automorphism of $H$. Then $H^\sigma$ is a hyperring. Moreover, if $H$ is a skew hyperfield, then $H^\sigma$ is a skew hyperfield. \end{lemma}
\proof We must first verify that $\boxplus$ is commutative and associative. Commutativity is clear from the symmetry in the definition. 
To see associativity, consider $aT^i, bT^j, cT^k$. If $i<j$, then 
$$(aT^i \boxplus bT^j)\boxplus cT^k= aT^i \boxplus cT^k = aT^i \boxplus (bT^j\boxplus cT^k)$$
If $i>j$ then 
$$(aT^i \boxplus bT^j)\boxplus cT^k= bT^j \boxplus cT^k = aT^i \boxplus (bT^j\boxplus cT^k)$$
So $i=j$, and by symmetry $j=m$. Then 
$$(aT^i \boxplus bT^i)\boxplus cT^i= (a+ b)T^i\boxplus cT^i=aT^i\boxplus (c+b)T^i= aT^i \boxplus (bT^i\boxplus cT^i).$$
Next, we show that $(H^\sigma, T^\infty,  \boxplus)$ satisfies the hypergroup axioms (H0), (H1).

(H0): $aT^i\boxplus T^\infty=  aT^i$ by definition.

(H1): For $aT^i\in H^\sigma$, we have $T^\infty\in aT^i\boxplus bT^j$ if and only if $i=j$ and $a=-b$. Thus $-(a\cdot T^i)=(-a)\cdot T^i$.

(H2): Suppose $a T^i\in bT^j\+cT^k$. We must show $cT^k\in a T^i\+ -bT^j$. If $j<k$, then $a T^i=bT^j$ and hence $cT^k\in a T^i\+ -bT^j$, and similar if $j>k$. If $i>j=k$, then $b=-c$ and hence $cT^k\in a T^i\+ -bT^j$.
So $i=j=k$ and $a\in b+c$, so that $c\in a+ (-b)$ and hence $cT^k\in a T^i\+ -bT^j$.

It is evident that $(H^\sigma\setminus\{T^\infty\}, T^0, \cdot)$ is a multiplicative monoid. We have $aT^i\cdot T^\infty=T^\infty$, so the zero $T^\infty$ is absorbing. Distributivity is straightforward.

Finally, if $H$ is a skew hyperfield, then $1\in H^\sigma$ is distinct from $0\in H^\sigma$, and each $aT^i$ has multiplicative inverse $\sigma^{-i}(a) T^{-i}$, so that $(H^\sigma\setminus\{T^\infty\}, T^0, \cdot)$ is a multiplicative group. Then $H^\sigma$ is a skew hyperfield.
\endproof
For any hyperring $H$, there is a hyperring homomorphism $\zeta: H^{\sigma}\rightarrow \Z_{\min}$ given by $\zeta: aT^i\mapsto i$ and $0\mapsto \infty$, and there is a group homomorphism $\tau: (\Z,+)\rightarrow (H^\sigma)^\star$ given by $\tau: i\mapsto T^i$. 

If we have $1\in 1+1$ in $H$, then $$\tau(\min\{i,j\})=T^{\min\{i,j\}}\in T^i\+ T^j=\min\{\tau(i), \tau(j)\}$$ even if $i=j$. For such $H$,  we may extend $\tau$ to a  hyperring homomorphism $\Z_{\min}\rightarrow H^\sigma$, by setting $\tau:\infty\mapsto 0$. 

If $H=\mathbb{K}$ and $\sigma$ is the identity, then $\zeta$ is an isomorphism with inverse $\tau$. In this sense, $\Z_{\min}\cong \mathbb{K}^{id}$ is a special case of the above construction.

\subsection{\label{ss:ore}Ore extensions of fields}
The definition of the above skew hyperfield of monomials was inspired by a construction of  skew fields due to Ore \cite{Ore1933}.

 Let $K$ be a skew field, and let $\sigma:K\rightarrow K$ be an automorphism. 
The {\em Ore extension} $K[T, \sigma]$ is the ring of formal polynomials $\sum_{i=0}^n a_i T^i$ in which $T$ commutes with elements $a\in K$ according to the rule $Ta=\sigma(a)T$. 

The ring $R=K[T, \sigma]$ satisfies left and right {\em Ore conditions}: for each $s,t\in R$, we have $sR\cap tR\neq \emptyset$ and $Rs\cap Rt\neq \emptyset$, which allows to define the
 left field of fractions $$K(T,\sigma):=\{a^{-1}b: a,b\in K[T,\sigma]\}.$$
 
There is  a hyperring homomorphism $\nu:K[T, \sigma]\rightarrow \mathbb{Z}_\infty$ determined by
$$\nu\left(\sum_{i=0}^n a_i T^i\right)=\min\{i: a_i\neq 0\}$$
and $\nu(0):=\infty$. 
This $\nu$ extends to $\nu:K(T,\sigma)\rightarrow \Z_{\min}$ by setting 
$\nu(a^{-1}b):=-\nu(a)+\nu(b)$.

If $K$ is a skew field, then there is a hyperring homomorphism $\mu: K[T,\sigma]\rightarrow K(T,\sigma,\min)$  determined by
$$\mu\left(\sum_{i=0}^n a_i T^i\right)=a_mT^m, \text{where }m=\min\{i: a_i\neq 0\}$$
This hyperring homomorphism extends to $\mu: K(T,\sigma)\rightarrow K(T,\sigma, \min)$ by setting  $\mu(a^{-1}b)=\mu(a)^{-1}\mu(b)$.
\begin{lemma} $\mu$ and $\zeta$ are hyperfield homomorphisms, and $\nu=\zeta\circ\mu$.\end{lemma}
There is a similar homomorphism $K(T,\sigma)\rightarrow K(T,\sigma,\max)$ which picks up the leading term.

\subsection{\label{ss:boundary}The boundary matroid of an $H^\sigma$-matroid}
Consider a $\Z_{\min}$-matroid $M$ on $E$ with Grassmann-Pl\"ucker function $\phi: E^r\rightarrow \Z_{\min}$. As the hyperaddition of $\Z_{\min}$ is idempotent, we have  $x=-x$ in $\Z_{\min}$ and hence the otherwise alternating Grassmann-Pl\"ucker function becomes oblivious to the ordering of its argument: $\phi(b_1,\ldots, b_r)=\phi(b_1',\ldots, b_r')$ whenever $\{b_1,\ldots, b_r\}=\{b_1',\ldots, b_r'\}$. Let $\nu:\binom{E}{r}\rightarrow \Z_{\min}$ be determined by
$$\nu(B):= \phi(b_1,\ldots, b_r)$$
whenever $B=\{b_1,\ldots, b_r\}$. Then $\nu$ is a {\em matroid valuation}, and it was shown by Dress and Wenzel \cite{DressWenzel1992a} that $$\mathcal{B}_0:=\left\{B\in \binom{E}{r}: \nu(B)=\min\left\{\nu(B'): B'\in \binom{E}{r}\right\}\right\}$$
is a nonempty set satisfying the base exchange axiom. We will call the matroid $M_0$ with ground set $E$ and set of  bases $\mathcal{B}_0$, the {\em boundary matroid of $M$}.\footnote{Dress and Wenzel speak of a {\em residue class geometry} in  \cite{DressWenzel1992a}.}

We will define boundary matroids more generally for  $H^\sigma$- matroids. Consider the natural hyperfield homomorphism $\zeta: H^\sigma\rightarrow \Z_{\min}$ given by $\zeta:aT^i\mapsto i$. 
 \begin{lemma} \label{lem:boundary}Let $H$ be a skew hyperfield and let $M$ be a left $H^\sigma$-matroid, and let $N:=(\zeta_* M)_0$. Let $[.]_0$  be the restriction of $[.]_M$ to $A_N$. Then $[.]_0$ are quasi-Pl\"ucker coordinates for $N$, taking values in $H$. \end{lemma}
 \proof Recall that by definition of the boundary matroid of a $\Z_{\min}$-matroid, the matroid $N$ has bases
 $$\mathcal{B}_0:=\left \{B\in \mathcal{B}: \nu(B)=\min\{\nu(B'): B'\in \mathcal{B}\}\right\},$$
where $\mathcal{B}$ is the set of bases of $\zeta_* M$ and $\nu(B)-\nu(B')=\zeta [B, B']$ for any $(B, B')\in A_M$.
Assuming without loss of generality that $\min\{\nu(B'): B'\in \mathcal{B}\}=0$, we have $\mathcal{B}_0:=\left \{B\in \mathcal{B}: \nu(B)=0\right\}$, and $\nu(B)>0$ if $B\in \mathcal{B}\setminus\mathcal{B}_0$. In particular $[B, B']_0=[B, B']_M\in H$ for all $(B,B')\in A_N$, since for such $(B,B')$ we have $\zeta([B, B']_M)=\nu(B)-\nu(B')=0$.

To prove that $[.]_0$ are quasi-Pl\"ucker coordinates for $N$, we need only show that $[.]_0$ satisfies (P3) by Lemma \ref{lem:weakimage}.
Consider $F, a,b,c,d$ so that $Fac, Fad, Fbc, Fbd\in \mathcal{B}_0$, and $Fab\not\in \mathcal{B}_0$. If $Fab$ is not a base of $M$, then 
 $$[Fac, Fbc]_0=[Fac, Fbc]_M=[Fad, Fbd]_M=[Fad, Fbd]_0,$$
 and likewise if $Fcd$ is not a basis of $M$. If on the other hand both $Fab, Fcd\in \mathcal{B}$, then $$1\in [Fac, Fad]_M\x [Fbd, Fbc]_M\+[Fcd, Fad]_M\x [Fab, Fbc]_M$$
 by the fact that the quasi-Pl\"ucker coordinates of $M$ satisfy (P4). As $Fac, Fad, Fbc, Fbd\in \mathcal{B}_0$, and $Fab, Fcd\not\in \mathcal{B}_0$, we have
 $$\zeta([Fac, Fad]_M\x [Fbd, Fbc]_M)=\nu(Fac)-\nu(Fad)+\nu(Fbd)-\nu(Fbc)=0$$ and
  $$\zeta([Fcd, Fad]_M\x [Fab, Fbc]_M)=\nu(Fcd)-\nu(Fad)+\nu(Fab)-\nu(Fbc)> 0.$$ 
 Then $1\in [Fac, Fad]_M\x [Fbd, Fbc]_M\+[Fcd, Fad]_M\x [Fab, Fbc]_M=\{[Fac, Fad]_M\x [Fbd, Fbc]_M\},$ and hence $[Fac, Fad]_0\x [Fbd, Fbc]_0=[Fac, Fad]_M\x [Fbd, Fbc]_M=1.$
 \endproof
 If $M$ is a left $H^\sigma$-matroid, then by Lemma \ref{lem:boundary} there exists a matroid $M_0$ such that $\underline{M_0}=(\zeta_* M)_0$ and $M_0$ is a weak image of $M$. Clearly, there can be at most one such matroid.
We will call $M_0$ the {\em boundary matroid of $M$}. By the Lemma, $M_0$ is a left $H$-matroid.

\newcommand{\F}{{\mathcal{F}}}


\section{\label{sec:algebraic}Matroids over hyperfields from algebraic matroids}
\subsection{Preliminaries on field extensions, algebraic matroids}
Let $K$ be a field, and $E$ be a finite set. We write $K[X_E]:=K[X_e:e\in E]$ for the polynomial ring over $K$ with a variable $X_e$ for each element of $E$, and $K(X_E)$ for its field of fractions. 
For a polynomial $q\in K[X_E]$, let $\underline{q}$ denote the smallest set $F$ so that $q\in K[X_F]$, i.e. $\underline{q}$ is the set of indices of variables which are mentioned in $q$. 
\begin{lemma} \label{lem:elimination}Let $I\subseteq K[X_E]$ be an ideal, and let $q,r\in I$ be irreducible over $K$. If $\underline{q}\neq \underline{r}$ and $e\in \underline{q}\cap \underline{r}$, then there exists a nonzero polynomial $s\in I$ such that $e\not\in \underline{s}\subseteq \underline{q}\cup \underline{r}$.\end{lemma}
If $L$ is an extension field of $K$, and $x_e\in L$ for $e\in E$, then $x_F$ is {\em algebraically dependent over $K$} if there is a nonzero polynomial $q\in K[X_F]$ so that $q(x)=0$ (when variables and values are both indexed by $E$, then $q(x)$ arises by substituting $X_e$ with $x_e$ for all $e\in E$). 
\begin{theorem} Let $L/K$ be a field extension, let $E$ be a finite set and let $x_e\in L$ for each $e\in E$. Let $\mathcal{C}$ be the set of inclusionwise minimal elements of 
$$\mathcal{A}:=\{F\subseteq E: x_F\text{ is algebraically dependent over }K\}\del \{\emptyset\}.$$
Then $(E,\mathcal{C})$ is a matroid.
\end{theorem}
\proof (MC0) and (MC1) hold for $\mathcal{C}$ as $\mathcal{C}$ is the set of inclusionwise minimal elements of an $\mathcal{A}\subseteq 2^E\setminus\{\emptyset\}$. 
We prove (MC2). Suppose $C, C'\in \mathcal{C}$ are distinct. Then there are polynomials $q,r\in K[X_E]$ so that 
$C=\underline{q}$ and $C'=\underline{r}$. If $q$ is reducible, some factor $q'$ of $q$ will have $\emptyset\neq \underline{q'}\subseteq\underline{q}$, and then $\underline{q'}=\underline{q}$ by minimality of $C=\underline{q}$ in $\mathcal{A}$. Replacing $q,r$ by such a factor if necessary, we may assume $q,r$ are irreducible over $K$. Consider the homomorphism $h: K[X_E]\rightarrow L$ which maps $h:X_e\mapsto x_e$, and let $I:=\ker(h)$. Then $q,r\in I$, and by Lemma \ref{lem:elimination}, there exists a polynomial $s\in I$ so that 
$$e\not\in \underline{s}\subseteq \underline{q}\cup \underline{r}=C\cup C'.$$
Then $A:=\underline{s}\in \mathcal{A}$, so that there is some $C''\in \mathcal{C}$ with $C''\subseteq A\subseteq C\cup C'\del\{e\}$, as required.
\endproof
We denote the matroid of the elements $x\in L^E$ in the field extension $L/K$ by $M(K, x)$. 
\begin{lemma}\label{lem:scaling} Let $L/K$ be a field extension, let $x\in L^E$, and let $h: K[X_E]\rightarrow L$ be the homomorphism which maps $h:X_e\mapsto x_e$. If $C$ is a circuit of 
$M(K, x)$ and $q\in K[X_C]$, then $K[X_C]\cap \ker(h) = q K[X_C]$ if and only if $q$ irreducible. Moreover, if  $q K[X_C]=q'K[X_C]$ then $q=\alpha q'$ for some $\alpha\in K^\star$.
\end{lemma}
We say that a polynomial $q$ as in the lemma {\em decorates} the circuit $C$ of $M(K,x)$. 

\subsection{The space of derivations} 
Let $R$ be any ring. A {\em derivation} of $R$ is a map $D:R\rightarrow R$ such that
\begin{itemize}
\item[(D0)] $D(1)=0$
\item[(D1)] $D(x+y)= D(x)+D(y)$ 
\item[(D2)] $D(x y)= D(x)y + xD(y)$
\end{itemize}
If $S\subseteq R$, then we say that a derivation $D$ is an {\em $S$-derivation}  if $D(s)=0$ for all $s\in S$.  

Consider a field extension $L/K$ and $x\in L^E$, and let $D$ be a $K$-derivation.
For any polynomial $q\in K[X_E]$ so that $q(x)=0$, we have  $D(q(x))=D(0)=0$. Applying (D0), (D1), (D2) to expand $D(q(x))$ we obtain
$\sum_{e\in E} \frac{\partial q}{\partial x_e} D(x_e)=0$.  
Here $\frac{\partial q}{\partial x_e}$ denotes the formal derivative $\frac{\partial q}{\partial X_e}$ as evaluated in $X_E=x$. It follows that \begin{equation}\label{eq:orth}d(q):= \left(\frac{\partial q}{\partial x_e}:e\in E\right)\perp \left(D(x_e): e\in E\right)=:D(x).\end{equation}
The following stronger statement is Theorem 5.1 of \cite{LangBook}. In the statement of this theorem, $q^D$ denotes the result of applying $D$ to each coefficient of $q\in K[X_E]$.
\begin{theorem}\label{thm:langder} Let $L/K$ be a field extension, let $x\in L^E$. Let $h:K[X_E]\rightarrow L$ be the homomorphism such that $h(X_e)=x_e$, and let $q_1,\ldots, q_t$ be a set of generators for $\ker(h)$. Suppose $D$ is a derivation of $K$. If $u\in L^E$ is such that for $i=1,\ldots, t$
$$0=q_i^D(x)+\sum_e  \frac{\partial q_i}{\partial x_e}u_e,$$
then there is one and only one derivation $D^*$ of $K(x_E)$ coinciding with $D$ on $K$, and such that $D^*(x_e)=u_e$ for every $e\in E$. 
\end{theorem} 
This theorem may be used to characterize 
$Der(K,x):=\{D(x):D \text{ a $K$-derivation of }K(x_E)\}.$
\begin{corollary} Let $L/K$ be a field extension, let $x\in L^E$. Then
$$Der(K,x)=\{d(q): q\text{ decorates a circuit of }M(K,x)\}^\perp.$$ 
\end{corollary}
\proof The polynomials decorating the circuits of $M(K,x)$ generate the kernel of $h$ as in the theorem. 
Apply the theorem to the trivial $K$-derivation $D$. Since $D$ is trivial, we have $q^D(x)=0$ for any decorating polynomial $q$. 
We obtain that $D^*$ is a $K$-derivation of $K(x_E)$ if and only if $D^*(x)\perp d(q)$ for each polynomial $q$ decorating a circuit of $M(K,x)$. 
\endproof
If $q\in K[X_E]$, then clearly $\underline{d(q)}\subseteq \underline{q}$, but equality need not hold if $K$ has positive characteristic $p$. We then have $e\in \underline{q}\del \underline{d(q)}$ if and only if $q$ can be written as a polynomial in $X_e^p$. The polynomial $q$ is {\em separable} in $X_e$ exactly if $e\in \underline{d(q)}$. 

If $k$ is any subfield of $L$ and $y\in L$, then $y$ is {\em separable} over $k$ if there is a polynomial $q\in k[Y]$ which is separable in $Y$ so that $q(y)=0$. The {\em separable closure} of $k$ in $L$ is $$k^{sep}:=\{y\in L: y\text{ separable over }k\}.$$
As a consequence of Theorem \ref{thm:langder}, any derivation of $k$ will extend uniquely to $k^{sep}$. 
\begin{corollary} Let $L/K$ be a field extension, let $x\in L^E$. Then $\dim Der(K,x)$ equals the rank of $M(K,x)$. 
\end{corollary}
\proof $B$ is a basis of $M(K,x)$ if and only $K(x_E)$ is algebraic over $K(x_B)$. Pick a basis $B$ so that the index $[K(x_E):K(x_B)^{sep}]$ is as small as possible. Then for each $e\in E\del B$, the circuit $C\subseteq B+e$ is decorated by a polynomial $q$ which is separable in $X_e$. If not, $q$ (being irreducible) is separable in some $f\in C-e\subseteq B$. Taking $B':=B+e-f$, we then have 
$K(x_B)^{sep}\subseteq K(x_{B'})^{sep}$, and the inclusion is strict since $x_e\not\in K(x_B)^{sep}$ and $x_e\in K(x_{B'})^{sep}$. Then 
$[K(x_E):K(x_{B'})^{sep}]<[K(x_E):K(x_B)^{sep}]$, contradicting the choice of $B$.

Consider values $u_e\in K(x_E)$ satisfying the condition of Theorem \ref{thm:langder}. 
Observe that upon fixing $u_f$ for each $f\in B$, the values of $u_e$ for $e\in E\del B$ are determined by the relation $0=\sum_e  \frac{\partial q}{\partial x_e}u_e,$ where $q$ is the polynomial decorating $C\subseteq B+e$, since $\frac{\partial q}{\partial x_e}\neq 0$. Hence $\dim Der(K,x)\leq |B|$. On the other hand the derivations $(D_e:=\partial/\partial x_e)_{e\in B}$ are independent, since  $D_e(x_f)\neq 0$ if and only if $e=f$, for all $e,f\in B$. It follows that $\dim Der(K,x)\geq |B|$ as well, and hence $\dim Der(K,x)= |B|=r(M(K,x))$.
\endproof
\ignore{
\begin{lemma} Let $L/K$ be a field extension, let $x\in L^E$, and let $M:=M(Der(K, x))$. Then $M$ is a weak image of $M(K, x)$, with $r(M)=r(M(K, x))$.
\end{lemma} 
\proof We must show (1) that each independent set of $M$ is independent in $M(K, x)$, and (2) that some basis of $M$ is a basis of $M(K, x)$.

To see (1), suppose $F$ is a dependent set of $M(K, x)$. Then $F$ contains some circuit $C$ of  $M(K, x)$, which is decorated by a polynomial $q$. As $q$ is irreducible, we have  $\underline{d(q)}\neq \emptyset$, and as  $d(q)\perp Der(K,x)$, $\underline{d(q)}$ is a dependent set of 
$M$. Then $F$ is dependent in $M$ as well, since
$\underline{d(q)}\subseteq \underline{q}= C\subseteq F$.

To see (2), let $B$ be a basis of $M$, and suppose $B$ is not a basis of $M(K, x)$. Then there is a hyperplane $H$ of $M(K, x)$ containing $B$, and it follows that the rank of $M(K(x_H), x_{E\del H})=M(K,x)/H$ is 1. Then there is at least one nonzero derivation $D\in Der(K(x_H), K(x_E))$. Then $D_x\in Der(K,x)$ but $\underline{D_x}\cap H=\emptyset$, contradicting that $B$ is a basis of $M$.   
\endproof
If $K$ and $L$ are fields of positive characteristic $p$, then it is possible that $M(Der(K, x))\neq M(K, x)$, for example if we choose $x=(y, y^p)\in (L\setminus K)^2$.
}

\subsection{The matroid of $\sigma$-derivatives}
Let $K\subseteq L$ be a field extension in positive characteristic $p$, let $E$ be a finite set, let $x\in L^E$, and put $N:=M(K, x)$. We will assume that $L$ is algebraically closed, and we write $\sigma: L\rightarrow L$ for the Frobenius automorphism $\sigma:x\mapsto x^p$. In what follows, we will create a left $L^\sigma$-signature for $N$ and a right $L^\sigma$-signature for $N^*$, aiming to showing orthogonality of these signatures. For brevity, we will not repeat our choice $E, K, L, \sigma$ in the lemmas of this section.

For a vector $u\in \mathbb{N}^E$, write $x^u=\prod_{e\in E} x_e^{u_e}$. Let $q=\sum_{u} q_u x^u\in K[X_E]$, and put $$m_e:=\max\{m\in\N: p^m\text{ divides }u_e\text{ for all }u\text{ such that }q_u\neq 0\}.$$
Then let $\overline{q}\in K[Z_E]$ be the polynomial such that  $q=\overline{q}\left(X_e^{p^{m_e}}: e\in E\right)$. The {\em $\sigma$-derivative} $d^\sigma(q):E \rightarrow L^\sigma$ is defined as
$$d^\sigma(q): e \mapsto \frac{\partial \overline{q}}{\partial z_e} T^{m_e}$$
where $z_e:=x_e^{p^{m_e}}$ for each $e\in E$.
Note that $\underline{d^\sigma(q)}=\underline{\overline{q}}=\underline{q}$, since by construction $\overline{q}$ is separable in each variable $Z_e$.  Let
$$\mathcal{C}_x:=\{\alpha\cdot d^\sigma(q): ~q\text{ decorates a circuit }C\text{ of }N, ~\alpha\in (L^\sigma)^\star\}.$$
\begin{lemma}\label{lem:Cx} Let $x\in L^E$. Then $\mathcal{C}_x$ is a left $L^\sigma$-signature of $N$.\end{lemma}
\proof We verify (C0), (C1), (C2) for $\mathcal{C}_x$. Clearly, (C0) and (C1) are true by construction. To see (C2), suppose $U, V\in\mathcal{C}_x$ are such that $\underline{U}\subseteq \underline{V}$. By definition of $\mathcal{C}_x$, we have $U=\alpha\cdot d^\sigma(q)$ and $V=\alpha'  \cdot d^\sigma(q')$ where $q$ decorates $C$ and $q'$ decorates $C'$, so that $\underline{U}=C,\underline{V}=C'$ both are circuits of $M(K, x)$, and hence $\underline{U}=\underline{V}$. It follows that $q$ and $q'$ both decorate the same circuit $C$ of $M(K,x)$. By Lemma \ref{lem:scaling}, there is a $\beta\in K^\star$ so that $q'=\beta\cdot q$. Then 
$$V = \alpha' \cdot d^\sigma(q') = \alpha'  \cdot\beta  \cdot d^\sigma(q) = \alpha' \cdot \beta \cdot \alpha^{-1} \cdot U,$$ 
as required. \endproof
On the dual side, for any $K$-derivation $D$ of $K(x_E)^{sep}$ we define $D^\sigma(x): E\rightarrow L^\sigma$ by setting
 $$D^\sigma(x): e\mapsto T^{m_e}D\left(x_e^{p^{-m_e}}\right),$$
where $m_e=\max\{m\in \N: x_e^{p^{-m_e}}\in K(x_E)^{sep}\}$. If $C$ is a cocircuit of $N$,  $H=E\del C$ is the complementary hyperplane, and $D$ is a nonzero $K(x_H)$-derivation $D$ of $K(x_E)^{sep}$, then $D(z)\neq 0$ for all $z\in K(x_E)$ such that $z^{p^{-1}}\not\in K(x_E)$. Hence  $\underline{D^\sigma(x)}=H$. We define
$$\mathcal{D}_x:=\left\{ D^\sigma(x) \cdot \beta: ~D\text{ a $K(x_H)$-derivation of $K(x_E)^{sep}$}, ~D\neq 0,  ~H\text{ hyperplane of } N, ~\beta\in (L^{\sigma})^\star\right\}$$
\begin{lemma} \label{lem:Dx}Let $x\in L^E$. Then $\mathcal{D}_x$ is a right $L^\sigma$-signature of $N^*$.\end{lemma}
\proof We verify (C0), (C1), (C2) for $\mathcal{D}_x$, noting that for a right signature we must reverse the order of multiplication in these axioms. As before, (C0) and (C1) are true by construction. We verify (C2). Let $U,V\in \mathcal{D}_x$ have $\underline{U}\subseteq \underline{V}$. Since both supports are cocircuits of $N$, we have  $\underline{U}=C= \underline{V}$ for some cocircuit $C$ of $N$, and with $H=E\del C$ there are nonzero $K(x_H)$-derivations $D, D'$ of $K(x_E)^{sep}$ and $\beta, \beta'\in L^\sigma$ so that $U=D^\sigma(x)\cdot \beta$ and $V=(D')^\sigma(x)\cdot \beta'$. Since the set of $K(x_H)$-derivations of $K(x_E)^{sep}$ is a vector space of dimension 1, there is an $\alpha\in K(x_E)^{sep}$ so that $D'=D\cdot \alpha$. Then 
$$ V=(D')^\sigma(x)\cdot \beta'=D^\sigma(x)\cdot \alpha\cdot \beta'=U\cdot \beta^{-1}\cdot \alpha\cdot \beta',$$
as required.
\endproof

\begin{lemma}\label{lem:CDscale}Let  $x, y\in L^E$ and $n\in\Z^E$ be such that $y_e=x^{p^{n_e}}$ for all $e\in E$, and let $\rho:E\rightarrow L^\sigma$ be given by $\rho: e\mapsto T^{n_e}$. Then 
$\mathcal{C}_x=\mathcal{C}_y^\rho$ and  $\mathcal{D}_y=\mathcal{D}_x^{\rho}$.
\end{lemma}
\proof It suffices to prove the lemma for $n=\one_{e_0}$, where $e_0\in E$ is some fixed element. So then $y_e=x^p_e$ if $e=e_0$ and $y_e=x_e$ otherwise. Denote $N:=M(K,x)=M(K,y)$.

Consider a circuit $C$ of $N$, and suppose $U\in\mathcal{C}_x$ has $\underline{U}=C$. Then $U=\alpha\cdot d^\sigma(q)$ for some $q\in K[X_C]$ decorating  $C$ and an $\alpha\in L^\sigma$. Let $m\in \N^E$ be such that $q=\overline{q}\left(X^{p^{m_e}}\right)$. There are two cases to consider. If $m_{e_0}>1$, then $q$ is a polynomial in $X_{e_0}^p$, and substituting $X_e$ with $Y_e^{1/p}$ if $e=e_0$ and $Y_e$ otherwise in $q$ gives a polynomial $q'\in K[Y_C]$. This polynomial $q'$ is irreducible, for any factorization of $q'$ would induce a factorization of $q$. Hence $q'$ decorates $C$ in $M(K,y)$. By construction of $q'$, we have $\overline{q'}=\overline{q}\in K[Z_E]$, and $q'=\overline{q}\left(Y^{p^{m'_e}}\right)$, where $m=m'+n$. Hence
$$U=\alpha\cdot d^\sigma(q)=\alpha\cdot \left(\frac{\partial \overline{q}}{\partial z_e}T^{m_e}\right)_e =  \alpha\cdot \left(\frac{\partial \overline{q}}{\partial z_e}T^{m'_e}\cdot \rho_e\right)_e=\alpha\cdot (d^\sigma(q'))^\rho\in \mathcal{C}_y^\rho.$$
If $m_{e_0}=0$, then construct $q'\in K[Y_C]$ from $q^p$ by substituting $X_e$ with $Y_e^{1/p}$ if $e=e_0$ and $Y_e$ otherwise. Again, any factorization of $q'$ would induce a factorization of $q$, and hence $q'$ decorates $C$. This time $\overline{q'}=\sigma(\overline{q})$, and $q'=\sigma(\overline{q})\left(Y^{p^{m'_e}}\right)$, where $m+\one_E=m'+n$. Hence
$$U=\alpha\cdot d^\sigma(q)=\alpha\cdot \left(T^{-1}\cdot \frac{\partial \sigma(\overline{q})}{\partial z_e}T^{m_e+1}\right)_e =  \alpha\cdot T^{-1}\cdot \left(\frac{\partial \overline{q'}}{\partial z_e}T^{m'_e}\cdot \rho_e\right)_e=\alpha\cdot T^{-1}\cdot (d^\sigma(q'))^\rho\in \mathcal{C}_y^\rho,$$
hence also $U\in \mathcal{C}_y^\rho$. It follows that  $\mathcal{C}_x\subseteq \mathcal{C}_y^\rho$, so that $\mathcal{C}_x=\mathcal{C}_y^\rho$ since both are $L^\sigma$-signatures of $N$.

Consider a hyperplane $H$ of $N$, and let $V\in \mathcal{D}_y$. Then there is a nonzero $K(y_H)$-derivation $D$ of $K(y_E)^{sep}$ and a $\beta\in L^\sigma$ so that $V=D^\sigma(y)\cdot \beta$. 
By definition of $D^\sigma(y)$, there is an $m\in \N^E$ so that $y^{p^{-m_e}}\in K(y_E)^{sep}$ and $D^\sigma(y_e)=T^{m_e}D\left(y^{p^{-m_e}}\right)$ for each $e\in E$. Again, there are two cases. If $m_{e_0}>0$, then 
$$x_{e_0}^{p^{-m_{e_0}+1}}=y_{e_0}^{p^{-m_{e_0}}}\in K(y_E)^{sep}.$$ 
Then $D$ is a derivation of $K(x_E)^{sep}$ as well, and $D^\sigma(x_e)=T^{m'_e}D\left(x^{p^{-m'_e}}\right)$ where $m'=m-n$. Hence
$$V=D^\sigma(y)\cdot \beta= \left(T^{m_e}D\left(y^{p^{-m_e}}\right)\right)_e\cdot \beta=\left(\rho_eT^{m'_e}D\left(x^{p^{-m'_e}}\right)\right)_e\cdot\beta=(D^\sigma(x)\cdot \beta)^\rho\in \mathcal{D}_x^{\rho}.$$
If $m_{e_0}=0$, then $x_e\not \in K(y_E)^{sep}$.  Then  $D':z\mapsto D(z^p)^{(1/p)}$ is a derivation of $K(x_E)^{sep}\subseteq \left(K(y_E)^{sep}\right)^{(1/p)}$, and taking $m'=m-n+\one_E$ we have 
$(D')^\sigma(x_e)=T^{m'_e}D'\left(x^{p^{-m'_e}}\right)$. Hence
 $$V=D^\sigma(y)\cdot \beta= \left(T^{m_e}D\left(y^{p^{-m_e}}\right)\right)_e\cdot \beta=\left(\rho_eT^{m'_e}D'\left(x^{p^{-m'_e}}\right)\cdot T^{-1}\right)_e\cdot\beta=((D')^\sigma_x\cdot T^{-1}\cdot \beta)^\rho\in \mathcal{D}_x^{\rho}.$$
and hence  $V\in \mathcal{D}_x^{\rho}$.
It follows that $\mathcal{D}_y\subseteq \mathcal{D}_x^{\rho}$ , so that $\mathcal{D}_y=\mathcal{D}_x^{\rho}$ since both are $L^\sigma$-signatures of $N^*$.\endproof

\begin{lemma} Let $x\in L^E$. Then $\mathcal{C}_x\perp\mathcal{D}_x$.\end{lemma} 
\proof Using Lemma \ref{lem:rescaling}, it is equivalent to prove that $\mathcal{C}_x^\rho\perp \mathcal{D}_x^{\rho^{-1}}$. We will invoke Lemma \ref{lem:CDscale} to simplify the argument. 

Let $U\in \mathcal{C}_x$ and $V\in \mathcal{D}_x$. 
Then $U=\alpha \cdot d^\sigma(q)$ for some circuit $C$ of $M(K,x)$ and $\alpha\in (L^\sigma)^\star$, and $V=D^\sigma(x)\cdot \beta$ for some   $K(x_H)$-derivation $D$, where $H$ is a hyperplane of $M(K,x)$ and  $\beta\in (L^{\sigma})^\star$. It is our object to prove that $U\perp V$, so that we may assume without loss of generality that $\alpha=\beta=1$.

By Lemma \ref{lem:CDscale}, we may assume that $V\in L^E$, and writing $U_e=T^{m_e}a_e$ with $a_e\in L$, that $$\min\{m_e: e\in \underline{U}\}=\min\{m_e: e\in \underline{U}\cap\underline{V}\}=0.$$
Then $D^\sigma(x)=D(x)$, and $d^\sigma(q)_e=d(q)_e$ for all $e\in \underline{U}$ so that $m_e=0$, so that
$$\sum_{e\in E} U_e V_e=\sum_{e\in \underline{U}\cap\underline{V}}  d^\sigma(q)_e \cdot D^\sigma(x_e)=
\sum_{e\in \underline{U}\cap\underline{V},m_e=0}  d^\sigma(q)_e\cdot  D^\sigma(x_e)=
 \sum_{e\in \underline{U}\cap\underline{V}}  d(q)_e \cdot D(x_e)=0,$$
as  the hypersum of any elements of $L^\sigma$ is determined by the terms $cT^m$ with $m$ minimal, and $d(q)\perp D(x)$.
\endproof

\begin{theorem} Let $K\subseteq L$ be a field extension in positive characteristic $p$, let $E$ be a finite set, let $x\in L^E$, and assume that  $L$ is algebraically closed.
Then $M:=(E, \mathcal{C}_x)$ is a left $L^\sigma$-matroid, and $M^*=(E, \mathcal{D}_x)$. 
\end{theorem}\proof Let $N:=M(K,x)$. By the lemma's of this section,  $\mathcal{C}_x$ is a left $L^\sigma$-signature of $N$, $\mathcal{D}_x$ is a right $L^\sigma$-signature of $N^*$,
and $\mathcal{C}\perp_3\mathcal{D}$. Using Theorem \ref{thm:duality}, it follows that 
$M:=(E, \mathcal{C}_x)$ is a left $L^\sigma$-matroid, and $M^*=(E, \mathcal{D}_x)$.
\endproof

We call the left $L^\sigma$-matroid $M^\sigma(K,x):=(E,\mathcal{C}_x)$ the {\em matroid of $\sigma$-derivatives}, and its dual $(E,\mathcal{D}_x)$ the {\em matroid of $\sigma$-differentials}, since each element $e$ of the ground set represents a differential $d(x_e)$. By construction, the matroid underlying $M^\sigma(K,x)$ is $M(K,x)$, but $M^\sigma(K,x)$ captures further information about $K, x$.

Recall the hyperfield homomophism $\zeta:L^\sigma\rightarrow \Z_{\min}$.
\begin{lemma} Let $K\subseteq L$ be a field extension characteristic $p>0$, let  $x\in L^E$ and assume that $L=K(x_E)$.
Let $M=M^\sigma(K,x)$. Then 
$$\zeta([Fa, Fb]_M)=log_p\frac{ [L: K(x_{Fa})^{sep}]}{[L: K(x_{Fb})^{sep}]}$$
for all bases $Fa, Fb$ of $M$.
\end{lemma}
\proof Let $q$ be the polynomial decorating the circuit $C$ so that $a,b\in C\subseteq Fab$. Suppose that 
$d^\sigma(q)_a=vT^i$ and $d^\sigma(q)_b=wT^j$. Then we have
$$p^{-i}[L: K(x_{Fb})^{sep}] = [L: K(x_{Fab})^{sep}] = p^{-j}[L: K(x_{Fa})^{sep}].$$
Also, $\zeta([Fa, Fb]_M)=-\zeta(vT^i)+\zeta(wT^j)=-i+j$. The lemma follows.\endproof
By a theorem of Cartwright \cite{Cartwright2017}, the Lindstr\"om valuation $\nu$ of $K,x$ is determined by 
$$\nu(B)=\log_p [K(x_E):K(x_B)^{sep}]$$
for each basis $B$ of $M(K,x)$. It follows that $\nu$ is a Grassmann-Pl\"ucker function for $\zeta_* M^\sigma(K,x)$.
\begin{lemma}Let $K\subseteq L$ be a field extension characteristic $p>0$, and let  $x\in L^E$.
Then $Der(K,x)$ is spanned by the cocircuits of $M^\sigma(K,x)_0$.
\end{lemma}
\subsection{Matroids over $K(T,\sigma)$}  If $K$ is a field of characteristic $p$ and $\sigma$ is the Frobenius map, then the elements the Ore ring $K[T,\sigma]$ naturally correspond to {\em $p$-polynomials}. Consider the map $\widehat{.}:K[T,\sigma]\rightarrow K[Z]$ given by
$$\widehat{\sum_{j} a_{j} T^j}= \sum_j a_j Z^{p^j}.$$
Then for any  $a,b\in K[T,\sigma]$, we have $\widehat{(a+b)}(Z)=\widehat{a}(Z)+\widehat{b}(Z)$ and $\widehat{ab}(Z)=\widehat{a}(\widehat{b}(Z))$.

For the remainder of this section, we consider a fixed field $K$, an extension field $L$ of $K$ and a transcendence base $z_1,\ldots, z_d$ of $L$ over $K$. In this context, there is a natural embedding 
$\psi:K(T, \sigma)^d\rightarrow L$, which sends vectors $v \in K(T, \sigma)^d$ to $p$-polynomials in $L$ as follows:
$$\psi: v\mapsto \sum_{i=1}^d\widehat{v_i}(z_i).$$
\begin{lemma}[Lindstr\"om\cite{Lindstrom1988}] Let $V\subseteq K[T,\sigma]^d$ be a finite set of vectors. Then $V$ is left linearly dependent over $K(T,\sigma)$ if and only if $\{\psi(v): v\in V\}$ is algebraically dependent over $K$.\end{lemma}
Let $E$ be a finite set and let $v_e\in  K[T,\sigma]^d$ for each $e\in E$. Let $M(v)$ be the left $K(T,\sigma)$-matroid which is linearly represented by the vectors $v_e$. With $x_e:=\psi(v_e)$ for all $e\in E$, we have $M(K, x)=\underline{M(v)}$ by Lindstr\"oms lemma. We show that in this context, the matroid of $\sigma$-derivatives $M^\sigma(K,x)$ may also be constructed directly from $M(v)$. Recall the skew field homomorphism $\mu: K(T,\sigma)\rightarrow K^\sigma$ from section \ref{ss:ore}, which maps $\mu: \sum_i a_iT^i\mapsto a_k T^i$, where $k=\min\{i: a_i\neq 0\}$. Let $\mu':K(T,\sigma)\rightarrow L^\sigma$ be given by $\mu'(a)=\mu(a)$.
\begin{lemma} Let $E\subseteq K[T,\sigma]^d$ be a finite set,  and let $x_e:=\psi(e)$ for all $e\in E$.  Then $M^\sigma(K, x)=\mu'_* M(v)$.\end{lemma}
\proof By Lindstr\"oms Lemma, we have $M(K,x)=\underline{M(v)}$, so that $M^\sigma(K, x)$ and $M(v)$ have the same underlying matroid. It therefore suffices to show that for each circuit $U$ of $M(v)$, the vector $\mu'_* U=(\mu'(U_e))_e$ is a circuit of $M^\sigma(K, x)$.

So consider a circuit  $U\in K(T,\sigma)^d$ of $M(v)$. By definition, $U$ is a left linear dependence $\sum_e U_ev_e=0$, of minimal support. 
Assume first that $U\in K[T,\sigma]^E$. Then the entries $U_e$ are formal polynomials in $T$, and we may define
$$q_U:=\sum_{e\in \underline{U}} \widehat{U_e}(X_e)\in K[X_E].$$
Since $U$ is a left linear dependence, we have $\left(\sum_e U_ev_e\right)_i=0$ for $i=1,\ldots, d$, and hence
$$q_U(x)=\sum_{e} \widehat{U_e}(x_e)= \sum_{e} \widehat{U_e}\left(\sum_i\widehat{(v_e)_i}\right) =\sum_e\sum_i \widehat{ (U_e v_e)_i}=\sum_i \widehat{ \left(\sum_e U_e v_e\right)_i} =0.$$
Hence, the polynomial $q_U$ decorates the circuit $\underline{U}$ of $M(K,x)$. We have $(d^\sigma q)_e = \mu'(\widehat{U_e})$ for each $e$, and hence $\mu'_* U=d^\sigma q$ is a circuit of  $M^\sigma(K,x)$.

In case $U\not\in K[T,\sigma]^E$, then there is a $c\in K[T,\sigma]$ so that $cU\in K[T,\sigma]^E$. Then $\mu_*(cU)$ is a circuit of $M^\sigma(K,x)$ and hence by the circuit axiom (C1), the vector $\mu'_*U=\mu'(c^{-1})\mu_*(cU)$ is a circuit of $M^\sigma(K,x)$.
\endproof

\subsection*{Example} Pick any $a\in K^\star$, and consider the following vectors from $K(T,\sigma)^2$: 
$$v_1=\left[\begin{array}{c} 1\\0\end{array}\right], v_2=\left[\begin{array}{c} 0\\T^3\end{array}\right], v_2=\left[\begin{array}{c}T^2+ T\\T^2\end{array}\right], v_4=\left[\begin{array}{c} 1\\T^4+aT\end{array}\right].$$ 
Each pair of these vectors is linearly independent over $K(T,\sigma)$, and hence $\underline{M(v)}\cong U_{2,4}$. Taking $x_e:=\psi(v_e)$ we obtain
$$x_1=z_1, ~x_2=z_2^{p^3}, ~x_3=z_1^{p^2}+ z_1^p+z_2^{p^2}, ~x_4=z_1+z_2^{p^4}+az_2^p.$$
where $L=K(z_1,z_2)$ has transcendence degree 2 over $K$. The circuit $U=(T^3+T^2, 1, -T, 0)$ gives rise to an algebraic relation $q_U=X_1^{p^3}+ X_1^{p^2}+X_2-X_3^p$, so that $d^\sigma(q_U)=(T^2, 1, -T, 0)=\mu_* U$. The  $K$-derivation $D=\frac{d}{dz_2}$ gives
$$D^\sigma(x)=(0, T^3, T^2, aT)=\mu_* V,$$
where $V=(0, T^3, T^2, T^4+aT)$ is a cocircuit of $M(v)$.

In $M^\sigma(K,x)$, we have the cross ratio $cr(1,2,3,4)=[13,14]\cdot[24,23]=\left(a^{-1}T\right)\cdot\left( T^{-1} \right)=a^{-1}\in L^\sigma$.

\section{\label{sec:flocks}Flocks}

\subsection{Preliminaries on matroid flocks}
In \cite{BDP2018}, Bollen, Draisma, and the present author defined a {\em matroid flock of rank $d$ on $E$} as a map $M$ which assigns a matroid $M_\alpha$ on $E$ of rank $d$ to each $\alpha\in \Z^E$, satisfying the following two axioms.
\begin{itemize}
\item[(MF1)] $M_\alpha/ i=M_{\alpha+\one_e}\del i$ for all $\alpha\in\Z^E$ and $e\in E$.
\item[(MF2)] $M_\alpha=M_{\alpha+\one_E}$ for all $\alpha\in \Z^E$.
\end{itemize}
Here, $\one_e$ denotes the unit vector in $\R^E$ with a 1 in the $e$-th position, and $\one_E$ the all-one vector in $\R^E$. More generally we write $\one_F:=\sum_{e\in F} \one_e$ for the incidence vector of any $F\subseteq E$.

Matroid flocks are cryptomorphic to {\em valuated matroids}. Using the definition of the boundary matroid from Section \ref{ss:boundary}, and noting that valuated matroids are essentially $\Z_{\min}$-matroids, we will now paraphrase this characterization, Theorem 7 of \cite{BDP2018}. Let $\mathbb{M}(E,r)$ denote the collection of matroids of rank $r$ on $E$.
\begin{theorem}\label{MF}Let $\mathcal{M}:\Z^E\rightarrow \mathbb{M}(E,r)$. The following are equivalent:
\begin{enumerate}
\item $\mathcal{M}$ is a matroid flock.
\item there is a $\Z_{\min}$-matroid $N$ so that $\mathcal{M}:\alpha\mapsto \left(N^{\alpha}\right)_0$.
\end{enumerate}
\end{theorem}
In what follows, we generalize this theorem to one that characterizes $H^\sigma$-matroids in terms of {\em $H^\sigma$-flocks}. In the proof, we will use one further lemma from \cite{BDP2018}.
If $M, M'$ are two matroids with common ground set $E$, then $M\succeq M'$ denotes that $M'$ is a {\em weak image} of $M$, i.e. that each dependent set of $M$ is also dependent in $M'$.  For any $\R_{\min}$-matroid $N$ on $E$, let 
$$C(N,\beta):=\left\{\alpha\in \R^E: \left(N^{\alpha}\right)_0\succeq \left(N^{\beta}\right)_0\right\}.$$
We will regard any $\Z_{\min}$-matroid as an $\R_{\min}$-matroid in the natural way.

The following is Lemma 18 of \cite{BDP2018}.
\begin{lemma}\label{lem:poly} Let $N$ be a $\Z_{\min}$-matroid on $E$ with valuation $\nu$, and let $\beta\in\Z^E$. Then $$C(N,\beta)=\left\{\alpha\in \R^E: \alpha_e-\alpha_f\geq \nu(B)-\nu(B-e+f)\text{ for all bases } B\text{ of } \left(N^{\beta}\right)_0, e\in B, f\in E\del B\right\}.$$
\end{lemma}
\subsection{$H^\sigma$-flocks and matroids over $H^\sigma$}
Let $H$ be a skew hyperfield, let $r\in \mathbb{N}$, and let $E$ be a finite set. Let $\mathbb{M}_H(E,r)$ denote the collection of left $H$-matroids of rank $r$ on $E$.
Consider an automorphism $\sigma$ of $H$.  An {\em $H^\sigma$-flock  of rank $r$ on $E$} is a map 
$\F:\Z^E\rightarrow \mathbb{M}_H(E,r)$, with the following properties:
\begin{itemize}
\item[(F1)] $\F_{\alpha+\one_e}\del e  = \F_\alpha/ e$ for all $\alpha\in \Z^E$ and $e\in E$.
\item[(F2)] $\F_{\alpha+\one_E} = \sigma_* \F_\alpha$ for all $\alpha\in \Z^E$.
\end{itemize}
We generalize Theorem \ref{MF}, which characterizes $\mathbb{K}^{id}$-flocks (matroid flocks) as cryptomorphic to $\mathbb{K}^{id}$-matroids ($\Z_{\min}$-matroids). In the proof of this§ generalization, we will use Theorem \ref{MF} itself as a stepping stone. Let $\tau: (\Z, +) \rightarrow (H^\sigma)^\star$ denote the group homomorphism $\tau: i\mapsto T^i$.
\begin{theorem} Let $\F:\Z^E\rightarrow \mathbb{M}_{H}(E,r)$. The following are equivalent:
\begin{enumerate}
\item $\F$ is an $H^\sigma$-flock.
\item there is a left $H^\sigma$-matroid $M$ so that $\F:\alpha\mapsto \left(M^{\tau(\alpha)}\right)_0$.
\end{enumerate}
\end{theorem}
\proof
(2)$\Rightarrow$(1): Assume (2). Let $N:=\zeta_* M$. Then $N$ is a $\Z_{\min}$-matroid, and therefore by Theorem \ref{MF},  the map 
$$\underline{\F}:\alpha\mapsto \left(N^{-\alpha}\right)_0=\underline{\left(M^{\tau(\alpha)}\right)_0}$$ is a matroid flock. We verify the two $H^\sigma$-flock axioms (F1) and (F2).

(F1): Without loss of generality, $\alpha=0$. We have $\F_0=M_0$, and $\F_{\one_e}=(M^\rho)_0$, where $\rho=\tau(\one_e)$.
As $\underline{\F}$ is a matroid flock, we have $$\underline{\F_0}\del e=\underline{M_0}\del e=\underline{(M^\rho)_0}/e=\underline{\F_{\one_e}}/e.$$
To show more strongly that $\F_0\del e=\F_{\one_e}/e$, it remains to show that also $[.]_{\F_0\del e}=[.]_{M_0\del e}=[.]_{(M^\rho)_0/e}=[.]_{\F_{\one_e}/e}$. 
If $e$ is not a coloop of  $M_0$, then for each  $(B, B')\in A_{M_0\del e}$ we have 
$$[B, B']_{M_0\del e}=[B, B']_{M}=[B, B']_{M^\rho}=[B, B']_{(M^\rho)_0/e}.$$
If $e$ is a coloop of  $M_0$, then $M_0\del e=M_0/e$, and for each  $(B, B')\in A_{M_0\del e}$ we have 
$$[B, B']_{M_0\del e}=[B+e, B'+e]_{M}=[B+e, B'+e]_{M^\rho}=[B, B']_{(M^\rho)_0/e}.$$
In either case, $[.]_{M_0}=[.]_{(M^\rho)_0}$, so that $\F_0=M_0=(M^\rho)_0=\F_{\one_e},$ as required.

(F2): Without loss of generality $\alpha=0$. Then $\F_0=M$, and $\F_{\one_E}=(M^\rho)_0$, where $\rho:=\tau(\one_E):e\mapsto T$.  
For each  $(B, B')\in A_{M_0}$, we have $$[B, B']_{M^\rho}=T[B, B']_M T^{-1}=\sigma([B, B']_M).$$ If $X$ is a circuit of $M$, then $T(X_eT^{-1}:e\in E)=(\sigma(X_e): e\in E)=\sigma(X)$ is a circuit of $M^\rho$. 
Hence 
$$\F_{\one_E}=(M^\rho)_0=\sigma_* (M_0)=\sigma_* \F_0,$$
as required.

(1)$\Rightarrow$ (2): Suppose (1). Then $\underline{\F}:\alpha\mapsto \underline{\F_\alpha}$ is a matroid flock. Hence by Theorem \ref{MF}, there is a $\Z_{\min}$-matroid $N$ so that
 $$\underline{\F_\alpha}= \left(N^{\alpha}\right)_0$$
If $M=(E,\mathcal{C})$ is a left $H^\sigma$-matroid so that $\F_\alpha=\left(M^{\tau(\alpha}\right)_0$, then the left quasi-Pl\"ucker coordinates $[.]=[.]_{\mathcal{C}}$ are a map $[.]:A_N\rightarrow H^\sigma$ so that $[.]_{\F_\alpha}=[.]^{\tau(\alpha)}$ for all $\alpha\in \Z^E$. That is, for each $\alpha\in\Z^E$
\begin{equation}\label{eq:flock}[Fa, Fb]_{\F_\alpha}= T^{\alpha_a}[Fa, Fb] T^{-\alpha_b}\end{equation}
whenever $Fa, Fb$ are adjacent bases of $\underline{\F_\alpha}$.  Conversely, if $[.]$ are left quasi-Pl\"ucker coordinates for $N$ satisfying these requirements, then $M:=(E,\mathcal{C}_{[.]})$ satisfies (2): then $\F_\alpha= \left(M^{\tau(-\alpha)}\right)_0$, as on either side of the equation, the matroids have the same underlying matroid and the same quasi-Pl\"ucker coordinates.

We first prove the existence of such a map $[.]:A_N\rightarrow H^\sigma$, satisfying \eqref{eq:flock} for each $\alpha$. So fix adjacent bases $Fa, Fb$ of $N$. We must argue that for each two $\alpha, \beta\in \Z^E$ so that $Fa, Fb$ are both bases of $\F_\alpha$ and $\F_\beta$, we have 
\begin{equation}\label{eq:consistent}T^{-\alpha_a}[Fa, Fb]_{\F_\alpha}T^{\alpha_b}=T^{-\beta_a}[Fa, Fb]_{\F_\beta}T^{\alpha_b}.\end{equation}
By (F2), we may assume that $\alpha\leq \beta$. We prove \eqref{eq:consistent} by induction on $\sum_e (\beta_e-\alpha_e)$. Let $e\in E$ be such that $\alpha_e<\beta_e$. If $e\in F$, then with $F':=F-e$, $F'a, F'b$ are  adjacent bases  of $\F_\alpha/e=\F_{\alpha+\one_e}\del e$, and hence
$$[Fa, Fb]_{\F_\alpha}=[F'a, F'b]_{\F_\alpha/e}=[F'a, F'b]_{\F_{\alpha+\one_e}\del e}=[Fa, Fb]_{\F_\beta}.$$
Taking $\alpha'=\alpha+\one_e$ and using the induction hypothesis on the pair $\alpha', \beta$, we obtain
$$T^{-\alpha_a}[Fa, Fb]_{\F_\alpha}T^{\alpha_b}=T^{-\alpha'_a}[Fa, Fb]_{\F_{\alpha'}}T^{\alpha'_b}=T^{-\beta_a}[Fa, Fb]_{\F_\beta}T^{\beta_b}.$$
If  $e\not\in Fab$, then $Fa, Fb$ are both bases of $\F_\beta\del e=\F_{\beta-\one_e}/e$, and hence
$$[Fa, Fb]_{\F_{\beta-\one_e}}=[Fa, Fb]_{\F_{\beta-\one_e}/ e}=[Fa, Fb]_{\F_\beta\del e}=[Fa, Fb]_{\F_\beta}.$$
Taking $\beta'=\beta-\one_e$ and again using induction, we have 
$$T^{-\alpha_a}[Fa, Fb]_{\F_\alpha}T^{\alpha_b}=T^{-\beta'_a}[Fa, Fb]_{\F_{\beta'}}T^{\beta'_b}=T^{-\beta_a}[Fa, Fb]_{\F_\beta}T^{\beta_b}.$$
Thus we have reduced to the case when $\alpha_e=\beta_e$ for all $e$ other than $a,b$. By Lemma \ref{lem:poly}, we have
$$C(N,\beta)=\left\{\alpha\in \R^E: \alpha_i-\alpha_j\geq \nu(B)-\nu(B-e+f)\text{ for all bases } B\text{ of } \left(N^{-\beta}\right)_0, e\in B, f\in E\del B\right\}.$$
Since $Fa, Fb$ are bases of of both $\F_\alpha$ and $\F_\beta$, it follows that $\alpha_a-\alpha_b\geq \nu(B)-\nu(B-e+f)=\beta_a-\beta_b$. Reversing $\alpha$ and $\beta$ in this argument, we also have $\beta_a-\beta_b\geq \nu(B)-\nu(B-e+f)=\alpha_a-\alpha_b$, so that $\alpha_a-\beta_a=\alpha_b-\beta_b$. It follows that  $\beta-\alpha=k\one_{ab}$. Consider the special case that $\beta-\alpha=\one_{ab}$, and let $G:=E\del ab$. We have
$$\F_\alpha/G=\F_{\alpha-\one_G}\del G=\F_{\alpha+\one_{ab}-\one_E}\del G=\sigma_*\F_{\alpha+\one_{ab}}\del G.$$
Then 
 $$T^{-\alpha_a}[Fa, Fb]_{\F_\alpha}T^{\alpha_b}=T^{-\alpha_a}[a, b]_{\F_\alpha/G}T^{\alpha_b}=T^{-\beta_a}[a, b]_{\F_{\alpha+\one_{ab}}\del G}T^{\beta_b}=T^{-\beta_a}[Fa, Fb]_{\F_\beta}T^{\beta_b}.$$
In general if $\beta-\alpha=k\one_{ab}$ with $k>1$, then  $\alpha':=\alpha+\one_{ab}\leq \beta$ and $\alpha'\in C(N, \beta)$, so that $Fa, Fb$ are bases of $\F_{\alpha'}$. The general case then follows by induction on $k$:
$$T^{-\alpha_a}[Fa, Fb]_{\F_\alpha}T^{\alpha_b}=T^{-\alpha'_a}[Fa, Fb]_{\F_{\alpha'}}T^{\alpha'_b}=T^{-\beta_a}[Fa, Fb]_{\F_\beta}T^{\beta_b}.$$
We have established that there exists a map $[.]:A_N\rightarrow H^\sigma$, satisfying \eqref{eq:flock} for each $\alpha$.

Next, we show that $[.]$ are left quasi-Pl\"ucker coordinates. Consider (P3), say. Suppose $Fac, Fad, Fbc, Fbd$ are bases of $N$, but $Fab$ or $Fcd$ are not. Then there exists an $\alpha\in \Z^E$ so that $Fac, Fad, Fbc, Fbd$ are bases of $\underline{\F}_\alpha$. By (P3) for $\F_\alpha$, we have 
$$[Fac, Fbc]=T^{-\alpha_a}[Fac, Fbc]_{\F_\alpha}T^{\alpha_b}=T^{-\alpha_a}[Fad, Fbd]_{\F_\alpha}T^{\alpha_b}=[Fad, Fbd].$$
To show (P0), (P1), (P2) it similarly suffices to that all bases in question are present in $\F_\alpha$ for some $\alpha\in\Z^E$.  

To show (P4), consider $F,a,b,c,d$ so that $\mathcal{B}':=\{Fac, Fad, Fbc, Fbd, Fab, Fcd\}$ are all bases of $N$. We need to show that $$1\in [Fbd, Fab]\x [Fac, Fcd]\+[Fad, Fab]\x [Fbc, Fcd].$$
Let $\nu$ be the valuation associated with  $N$, so $\zeta([B,B'])=\nu(B)-\nu(B')$ for all adjacent bases $B, B'$ of $N$. By Theorem \ref{thm:gp}, 
we have $\infty\in (\nu(Fab)+\nu(Fcd))\+(\nu(Fac)+\nu(Fbd))\+(\nu(Fad)+\nu(Fbc))$ in $\Z_{\min}$. That is, the minimum of the three numbers 
$$\nu(Fab)+\nu(Fcd), \nu(Fac)+\nu(Fbd), \nu(Fad)+\nu(Fbc)$$
is attained at least 2 times. There are four cases to consider. If   $\nu(Fab)+\nu(Fcd)= \nu(Fac)+\nu(Fbd)= \nu(Fad)+\nu(Fbc)$, then there exists an $\alpha$ so that $\mathcal{B}'\subseteq \F_\alpha$, and then (P4) holds as it holds in $\F_\alpha$. If $\nu(Fab)+\nu(Fcd)= \nu(Fac)+\nu(Fbd)< \nu(Fad)+\nu(Fbc)$, then 
$[Fbd, Fab]\x [Fac, Fcd]=1$ as there exists an $\alpha\in \Z^E$ so that $Fbd, Fab, Fac, Fcd$ are bases of $\F_\alpha$, and $Fad$ or $Fbc$ are not. Also,
$$\zeta([Fbd, Fab]\x [Fac, Fcd])=\nu(Fbd)-\nu(Fab)+\nu(Fac)-\nu(Fcd)=0$$
and $$\zeta([Fad, Fab]\x [Fbc, Fcd])=\nu(Fad)-\nu(Fab)+\nu(Fbc)-\nu(Fcd)>0$$
so that $$1=[Fbd, Fab]\x [Fac, Fcd] ~\in ~[Fbd, Fab]\x [Fac, Fcd]\+[Fad, Fab]\x [Fbc, Fcd]. $$ 
The case when $\nu(Fab)+\nu(Fcd)=\nu(Fad)+\nu(Fbc)< \nu(Fac)+\nu(Fbd)$ is similar.

If $\nu(Fab)+\nu(Fcd)> \nu(Fac)+\nu(bd)= \nu(ad)+\nu(bc)$, then $[Fac, Fad]= [Fbc, Fbd]$ and $[Fac,Fbc]=[Fad, Fbd]$ as before. Then $[Fbd, Fab]\x [Fac, Fcd]= -[Fad, Fab]\x [Fbc, Fcd]$, and 
$\zeta([Fbd, Fab]\x [Fac, Fcd])=\zeta([Fad, Fab]\x [Fbc, Fcd])<0$, so that 
$$1\in [Fbd, Fab]\x [Fac, Fcd]\+[Fad, Fab]\x [Fbc, Fcd],$$
as required.
\endproof

\section{Final remarks}
We showed that the theory of matroids over hyperfields can be extended to skew hyperfields, and that a algebraic matroid representation $K, x$ gives rise to a matroid $M^\sigma(K,x)$ over a skew hyperfield. This matroid comprises the same information as the Frobenius flock of \cite{BDP2018}, but the quasi-Pl\"ucker coordinates of $M^\sigma(K,x)$ give a better perspective of the overall coherence (and finiteness) of the data presented by such flocks. We hope that this will shed light on the hard problem of characterizing the algebraic representability of matroids. 

Many issues came to mind when writing this paper, which were outside the main scope.  We like to end by listing some of them.

\subsection{Cross ratios}

\ignore{Several ways present themselves for extending matroids over skew hyperfields as developed in Section \ref{sec:skew}. 

We conjecture that there is no real obstacle to extending the equivalence of circuit axioms, quasi-Plucker axioms, and dual signatures, to infinite matroid in the sense of \cite{Diestel2013}. The algebraic conditions our axioms impose are inherently local to minors of rank and corank at most 2, in the sense of Lemma \ref{lem:minor}.
Under the assumption that the underlying object is an (infinite) matroid, it should be possible to show the equivalences between these axiom systems at the level of these finite minors. 

Defining a Grassmann-Plucker function for an infinite matroid over a commutative hyperfield is perhaps more involved, but seems possible also. Considering the equivalence of bases $B\sim B:\Leftrightarrow |B\del B'|<\infty$, it should be possible to fix such a function after choosing the value of one representative of each equivalence class.

Our matroids over skew hyperfields correspond to what Baker and Bowler call {\em weak} matroids over hyperfields. The stronger notion is characterized by the duality condition $\mathcal{C}\perp \mathcal{D}$ (instead of $\mathcal{C}\perp_3 \mathcal{D}$), where $\mathcal{C}$ and $\mathcal{D}$ are the circuits and cocircuits of a matroid over a skew hyperfield.  We expect that the equivalent strong version of (C3) is  the same as in the commutative case, and that the corresponding strong version of (P4) is the obvious one. For infinite matroids, these stronger duality conditions seem very difficult to state even, unless the hyperfield in question can gracefully handle taking the hypersum of infinitely many terms.

}

In Section \ref{ss:cr}, we defined cross ratios and listed several of their properties. It is not clear to what extent these properties define matroids over skew hyperfields. There may not be coordinates which correspond with given cross ratios even if $N=U_{2,4}$ and $H$ is commutative. We conjecture that any obstacles to finding corresponding quasi-Pl\"ucker coordinates will be local, in the following precise sense. 

Let $N$ be a matroid. We say that a map $cr:CR_N\rightarrow H$ is {\em consistent} if there exists  quasi-Pl\"ucker coordinates $[.]$ for $N$ such that $cr(F,a,b,c,d)=[Fac,Fad][Fbd, Fbc]$ for all $(F,a,b,c,d)\in CR_N$. 
If $S$ is an independent set of $N$, and $T$ is an independent set of $N^*$, then a map $cr:CR_N\rightarrow H$ induces a map $cr/S\del T:CR_{N'}\rightarrow H$ on the minor $N':=N/S\del T$, determined by  $cr/S\del T: (F,a,b,c,d)\mapsto cr(S\cup F,a,b,c,d)$.
\begin{conjecture}Let $N$ be a matroid on $E$, and let $H$ be a skew hyperfield such that $1=-1$ if $N$ has a Fano minor. The following are equivalent for any map $cr:CR_N\rightarrow H$:
\begin{enumerate}
\item $cr$ is consistent; and
\item $cr/S\del T$  is consistent for each disjoint $S,T\subseteq E$ so that $S$ is an independent set of $N$, $T$ is an independent set of $N^*$, and $N/S\del T$ has at most 5 elements.  
\end{enumerate}
\end{conjecture}
The special case of this conjecture where $H=\mathbb{S}$ is a theorem of Gelfand, Rybnikov, and Stone \cite{GRS1995}, and if $H$ is commutative the conjecture follows from the work of Delucchi, Hoessly, Saini \cite{DHS2018}.

\subsection{The skew hyperfield of monomials}
 If $H$ is a field and $\sigma$ is the identity, then $H(T,\sigma,\max)$ is commutative and equals the {\em hyperfield of monomials} described by Viro in \cite{Viro2010}.
Viro notes that the role of $\Z$ in his definition can be replaced by any linearly ordered group $(\Gamma,+,<)$. This seems to apply also to our construction. Consider a skew hyperfield $H$, and automorphism $\sigma_i$ of $H$ for each $i\in \Gamma$ so that $\sigma_{i+ j}=\sigma_i\circ\sigma_j$ for all $i,j\in\Gamma$. Then we can define a hyperfield 
$$H\rtimes_\sigma \Gamma_{\max}:=(H\times \Gamma, 1,0,\odot,\+)$$
with  $1:=(1,0)$, $0:=(0,0)=(0,i)$ for all $i\in \Gamma$, multiplication given by  
$(a,i)\odot(b,j):=(a\sigma_i(b),i+ j)$
for all $a,b\in H$ and $i,j\in \Gamma$,  and addition given by $0\+x= x\+ 0=\{x\}$ and 
$$(a,i)\boxplus (b,j):=\left\{\begin{array}{ll} 
\{(a,i)\}&\text{if } i>j\\
\{(b,j)\}&\text{if } i<j\\
(a+b)\times\{i\}&\text{if } i=j\text{ and } a\neq -b\\
(a+b)\times\{i\}\cup H^\star\times\{ k\in \Gamma, k< i\}& \text{if }i=j\text{ and }a=-b
\end{array}
\right. 
$$
for $a,b\in H^\star$ and $i,j\in \Gamma$, where $+$ is the hyperaddition of $H$. There is an obvious variant $H\rtimes_\sigma \Gamma_{\min}$.

This skew hyperfield resembles the {\em extended tropical hyperring} of \cite{AGG2014}, but it is different when adding  $(a,i)\boxplus (b,j)$ in the case that $i=j\text{ and } a\neq -b$. With trivial automorphisms $\sigma_i=id$ we have $\mathcal{T}(\R) \cong \mathbb{S}\rtimes \R_{\max}$ and $\mathcal{T}(\mathbb{C})\cong \Phi \rtimes \R_{\max}$. Here $\mathcal{T}(\mathbb{R})$ and $\mathcal{T}(\mathbb{C})$ are Viro's tropical reals and tropical complex numbers, and $\Phi$ is the tropical phase hyperfield.

For any $H, \Gamma$, there is a homomorphism $\zeta:H\rtimes_\sigma \Gamma_{\min}\rightarrow \Gamma_{\min}$ given by $\zeta: (a,i)\mapsto i$. 
 \begin{lemma} \label{lem:boundary2}Let $M$ be a left $H\rtimes_\sigma \Gamma_{\min}$-matroid, and let $N:=(\zeta_* M)_0$. Let $[.]_0$  be the restriction of $[.]_M$ to $A_N$. Then $[.]_0$ are quasi-Pl\"ucker coordinates for $N$, taking values in $H$. \end{lemma}
The proof of Lemma  \ref{lem:boundary}, which is the special case  of this statement where $\Gamma=\Z$, will also serve as a proof of Lemma \ref{lem:boundary2}.
Thus we may define the boundary matroid $M_0$ of any left $H\rtimes_\sigma\Gamma_{\min}$-matroid $M$ as the unique left $H$-matroid $M_0$ so that $\underline{M_0}=(\zeta_* M)_0$, and so that $M_0$ is a weak image of $M$.

In general, there is a hyperfield homomorphism $\upsilon: H\rightarrow H\rtimes_\sigma \Gamma_{\min}$ given by $\upsilon: h\mapsto (h,0)$ and a group homomorphism $\tau: \Gamma_{\min}\rightarrow H\rtimes_\sigma \Gamma_{\min}$ given by $\tau: i\mapsto (1,i)$. 
Thus any left $H$-matroid $M$ can also be considered as a left $H\rtimes_\sigma \Gamma_{\min}$-matroid $\upsilon_* M$. By rescaling and taking the boundary matroid, this matroid spawns a left $H$-matroid $\left((\upsilon_*M)^{\tau(w)}\right)_0$ for each $w\in \Gamma^E$, as in the characterization of flocks.  
This construction is not without precedent. Ardila, Klivans and Williams \cite{AKW2006} define, for any oriented matroid $M$ on $E$ and $w\in\R^E$, an {\em initial oriented matroid} $M_w$ on $E$. We have $M_w=\left((\upsilon_*M)^{\tau(w)}\right)_0$, where $\tau$ and $\upsilon$ are the canonical maps for $H=\mathbb{S}$ and $\Gamma=\R$. That $M_w$ is indeed an oriented matroid then follows from Lemma \ref{lem:boundary2}.

\subsection{Groebner bases in positive characteristic} In Section \ref{sec:algebraic}, we considered a field $K$ of positive characteristic $p$, an extension field $L$ and elements $x_e\in L$ for $e\in E$. The results of this section highlight that $K(x_E)$ has a certain robustness against applications of the Frobenius map $\sigma: x\mapsto x^p$ to the individual elements $x_e$. If $y_e = x_e^{p^{m_e}}$, then for any irreducible $q\in K[X_E]$ so that $q(x)=0$ there is an irreducible $q'\in K[Y_e]$ so that $q'(y)=0$, and $q^{p^n}(X_e:e\in E)=q'(X_e^{p^{m_e}}:e\in E)$. That is, irrespective of such Frobenius actions, the irreducible polynomial relations are always just a variation of the same polynomial $\overline{q}\in K[Z_E]$.

In the light of this invariance, it seems inappropriate that of a Groebner basis would change more than superficially when substituting a variable $X$ by $X^p$, or that the steps taken by the Buchberger algorithm would turn out truly different. We imagine a variant which is indifferent to such changes.

To make the Buchberger algorithm ignore substitutions such as the above, we may no longer distinguish between a polynomial $q$ and its power $q^{p^k}$. The monomial order $\preceq$ on $\N^E$ must ignore powers of $p$. That is, for any $u, v\in \N^E$ we must have 
$$u\preceq v\text{ if and only if }u'\preceq v'$$
where $u_e=u'_e p^{val_p(u_e)}$, $v_e=v' p^{val_p(v_e)}$ for each $e\in E$. The monomial order could otherwise be lexicographic, based on a linear order $<$ of $E$.
When using $q$ with leading monomial $X^u$ to reduce $r$ with leading monomial $X^v$, we must first replace $q$ with a $p^k$-th power to ensure that $val_p(u_e)=val_p(v_e)$, where $e=\max\{f\in E: u_e\neq 0, v_e\neq 0\}$. Here the maximum is taken with respect to the chosen order of $E$.

We are not aware of any such variant of the Buchberger algorithm in the literature, but we think this could be the more efficient way to decide independence of sets in algebraic matroids. 

\ignore{

(TODO): relation with Evans-Hrushovski lemma on $M(K_4)$, and with Lindstr\"oms harmonic conjugate.
\begin{conjecture} $cr(a,b,c,d) = 1$ in an algebraic matroid $\Leftrightarrow$ the representation is equivalent to a one derived from an algebraic group. 
\end{conjecture}

}
\bibliographystyle{alpha}
\bibliography{math}

\end{document}